\documentclass[12pt]{amsart}
\title[Relative Schottky sets]{
Planar relative Schottky sets and quasisymmetric maps
}%
\author{Sergei Merenkov\\ \\ \textit{Dedicated to the memory of Oded Schramm}}
\address{Department of Mathematics\\
University of Illinois\\ 1409 W Green St\\ Urbana, IL
61801\\USA} \email{merenkov@illinois.edu}
\thanks{Supported by NSF grants DMS-1001144, DMS-0703617, DMS-0653439, DMS-0400636.}

\date{\today}

\usepackage[active]{srcltx}

\usepackage{amsmath, amssymb, amsthm,latexsym}
\usepackage{graphicx}
\newcommand\C{{\mathbb C}}

\newcommand\N{{\mathbb N}}

\newcommand\R{{\mathbb R}}

\newcommand\Sph{{\mathbb S}}
\newcommand\dee{\partial}

\newcommand\id{\operatorname{id}}

\newcommand\co{\colon}
\renewcommand\:{\colon}
\newcommand\sub {\subseteq}
\newcommand\ra {\rightarrow}

\def\area{\mathop{\mathrm{area}}}

\newcommand\no{\noindent}

\newtheorem{theorem}{Theorem}[section]

\newtheorem{conjecture}[theorem]{Conjecture}
\newtheorem{proposition}[theorem]{Proposition}
\newtheorem{corollary}[theorem]{Corollary}

\newtheorem{remark}{Remark}
\newtheorem{lemma}[theorem]{Lemma}

\theoremstyle{definition}


\begin{document}

\abstract
{A relative Schottky set in a planar domain $\Omega$ 
is a subset of $\Omega$ obtained by removing from $\Omega$ open geometric discs whose closures are in $\Omega$ and are pairwise disjoint. 
In this paper we study quasisymmetric and related maps between relative Schottky sets of measure zero. 
We prove, in particular, that quasisymmetric maps between such sets in Jordan domains are conformal,
locally bi-Lipschitz, and that their first derivatives are locally Lipschitz. 
We also provide a locally bi-Lipschitz uniformization result for relative Schottky sets 
in Jordan domains and establish rigidity with respect to local quasisymmetric maps for relative Schottky sets in the unit disc. 
}
\endabstract

\maketitle

\section{Introduction}\label{s:Intro}

\no
Let $\Omega$ be a domain
in the standard $n$-sphere 
$$
\Sph^n=\{(x_1,x_2,\dots, x_{n+1})\in\R^{n+1}\: |x_1|^2+|x_2|^2+\dots+|x_{n+1}|^2=1\}. 
$$
A
\emph{relative Schottky set} $S$ in $\Omega$ is a subset of $\Omega$ whose
complement in $\Omega$ is a union of open geometric  balls $\{B_i\}_{i\in I}$ with closures $\overline B_i,\ i\in I$, in $\Omega$, and such that $\overline B_i\bigcap\overline B_j=\emptyset,\ i\neq j$. 
We write
$$
S=\Omega\setminus \bigcup_{i\in I}B_i.
$$
The
boundaries of the balls $B_i$ are called \emph{peripheral
spheres} or, if $n=2$, \emph{peripheral circles}. 
If $\Omega$ is the sphere $\Sph^n$ or the Euclidean space $\R^n$, a relative
Schottky set in $\Omega$ is called a \emph{Schottky set}. Schottky sets arise in geometry as boundaries at infinity of universal covers of compact hyperbolic manifolds with non-empty totally geodesic boundaries. 
Relative Schottky sets, endowed with the restriction of the spherical
metric, were introduced in~\cite{BKM07} in connection with quasisymmetric rigidity.  
The main purpose of this paper is to investigate local and infinitesimal 
properties of quasisymmetric maps between relative Schottky sets in 
Jordan domains contained in $\Sph^2$.

Let $(X,d_X)$ and $(\tilde X,d_{\tilde X})$ be metric spaces and $\eta\:
[0,\infty)\to[0,\infty)$ be an arbitrary homeomorphism. A
homeomorphism $f\: X\to \tilde X$ is called $\eta$-\emph{quasisymmetric}
if
$$
\frac{d_{\tilde X}(f(p),f(q))}{d_{\tilde X}(f(p),f(r))}\leq\eta
\bigg(\frac{d_X(p,q)}{d_X(p,r)}\bigg),
$$
for every triple of distinct points $p,q$, and $r$ in $X$. A
homeomorphism between metric spaces is called
\emph{quasisymmetric} if it is $\eta$-quasisymmetric for some
$\eta$. 
We say that a homeomorphism between two metric spaces $X$ and $\tilde X$ is \emph{locally quasisymmetric} if its restriction to every compact set $K$ in $X$ is $\eta_K$-quasisymmetric, with $\eta_K$ depending on $K$.

A \emph{M\"obius transformation} in $\Sph^n$
is a composition of finitely many reflections in $(n-1)$-spheres in
$\Sph^n$. 
The image of every relative Schottky set under a M\"obius transformation is a relative Schottky set.
Every M\"obius transformation is 
quasisymmetric.
Let $S$ be a relative Schottky set and $\mathcal F$ be a family of deformations of $S$, so that every $f\in\mathcal F$ is a homeomorphism of $S$ onto a relative Schottky set $\tilde S$ that may depend on $f$. E.g., $\mathcal F$ may consist of all  quasisymmetric or locally quasisymmetric deformations.
A relative Schottky set $S$ is called \emph{rigid} with respect to $\mathcal F$ if every $f\in \mathcal F$ is
the restriction to $S$ of a M\"obius transformation. 
%
%
%
%
The following three theorems were proved in~\cite{BKM07}. 
\smallskip\no
\newline
{\bf Theorem A.} 
Every Schottky set in $\Sph^n,\ n\in\N, n\geq 2$, of spherical measure zero is rigid with respect to  quasisymmetric maps.
\smallskip\no
\newline
{\bf Theorem B.}
A Schottky set in $\Sph^2$ is rigid with respect to  quasisymmetric maps if and only if it has spherical measure 
zero.
\smallskip\no
\newline
{\bf Theorem C.} Let $n\in\N, n\geq 3$, and $\Omega\subseteq\Sph^n$.  Then every locally porous relative Schottky set in $\Omega$ is rigid with respect to quasisymmetric maps.

Roughly speaking, local porosity means that locally the peripheral spheres appear on all scales and locations. See~\cite{BKM07} for the definition.
The proof of Theorem~C shows that locally porous relative Schottky sets in domains in $\Sph^n,\ n\geq3$, are rigid with respect to locally quasisymmetric maps.
In contrast, the following theorem, generalizing Theorem~B, shows that rigid with respect to quasisymmetric maps relative Schottky sets  in
domains contained in $\Sph^2$ form a narrow class.
\begin{theorem}\label{T:Crit}\
A relative Schottky set $S$ in $\Omega\subseteq\Sph^2$ is rigid with respect to quasisymmetric maps if and only if $S\bigcup(\Sph^2\setminus\Omega)$ has spherical measure zero.
\end{theorem}

Recall that if $(X,d_X)$ and $(\tilde X, d_{\tilde X})$ are metric spaces, a map $f\: X\to \tilde X$ is said to be 
$L$-\emph{Lipschitz}, $L>0$, if
$$
d_{\tilde X}(f(p),f(q))\leq Ld_X(p,q),
$$
for all $p,q\in X$. We say that $f\:X\to\tilde X$ is \emph{locally Lipschitz} if every $p\in X$ has a neighborhood $U$ and a constant $L$ such that $f$ restricted to $U$ is $L$-Lipschitz. 
A homeomorphism $f\: X\to\tilde X$ is called $L$-\emph{bi-Lipschitz},
$L\geq 1$, if 
$$
\frac1{L}d_{X}(p,q)\leq d_{\tilde X}(f(p),f(q))\leq Ld_{X}(p,q),
$$ 
for all $p,q\in X$. We say that a homeomorphism $f\:X\to\tilde X$ is \emph{locally bi-Lipschitz}, if  every $p\in X$ has a neighborhood $U$ and a constant $L$ such that $f$ restricted to $U$ is an $L$-bi-Lipschitz homeomorphism onto its image.

The following theorem 
establishes conformality and the local bi-Lipschitz property of quasisymmetric maps between relative Schottky sets in Jordan domains in the plane.
\begin{theorem}\label{T:Diff}
Suppose that $S$ is a relative Schottky set of measure zero in a Jordan domain 
$\Omega\subseteq\C$. 
Let $f\: S\to \tilde S$ be a locally quasisymmetric orientation preserving map from $S$ to a relative Schottky set $\tilde S$ in a Jordan
domain $\tilde\Omega\subseteq\C$. Then 
$f$ is conformal in $S$ in the sense that for every 
$p\in S$,
\begin{equation}\label{E:Der}
f'(p)=\lim_{q\to p,\, q\in S}\frac{f(q)-f(p)}{q-p}
\end{equation}
exists and is not equal to zero.
Moreover, the map $f$ is locally bi-Lipschitz in $S$ and the first derivative of $f$ defined by~(\ref{E:Der}) is continuous in $S$. 
\end{theorem}

If $S$ is locally porous, it is not hard to see using standard compactness arguments 
that for every $p\in S$ there are two sequences of scales $(r_k)$ and $(\tilde r_k)$, $0<r_k, \tilde r_k\to0$ as $k\to\infty$, with the following properties.
The sequences of sets 
$
(S_k={(S-p)}/{r_k})_{k\in\N}
$
and
$(\tilde S_k={(\tilde S-f(p))}/{\tilde r_k})_{k\in\N}
$
converge in the Gromov-Hausdorff topology to Schottky sets $S_p$ and $\tilde S_{f(p)}$, 
respectively, and  the sequence of maps 
$
(q\mapsto{(f(p+r_kq)-f(p))}/{\tilde r_k})_{k\in\N}
$
from $S_k$ to $\tilde S_k$ converges locally uniformly to a quasisymmetric map $f_p$ from $S_p$ to $\tilde S_{f(p)}$. 
An application of Theorem~\ref{T:Crit} shows that $f_p$ is the restriction of a conformal linear map. 
The conclusion of Theorem~\ref{T:Diff} is much stronger in the sense that the limit in~(\ref{E:Der}) is independent of sequences of scales. 

We believe that quasisymmetric maps between relative Schottky sets of measure zero possess higher degree of regularity as in the following conjecture, motivated by~\cite{HS98}.
\begin{conjecture}\label{C:Diff}
Let $f\: S\to \tilde S$ be an orientation preserving  quasisymmetric map between relative Schottky sets (not necessarily in Jordan domains) of measure zero. Then $f$ is conformal at each point $p\in S$ and $f\in C^{\infty}(S)$, i.e., the derivatives of $f$ of all orders exist on $S$ in the sense of equation~(\ref{E:Der}).
\end{conjecture}

Under slightly stronger assumptions than those of Theorem~\ref{T:Diff}, we can prove the following quantitative statement akin to the Koebe distortion theorem for  Riemann maps between Jordan domains.
If $X$ is a metric space, $p\in X$, and $r>0$, let $B(p,r)$ denote the open ball in $X$ of radius $r$ centered at $p$. 
\begin{theorem}\label{T:Quantlip}
If $f$ as in Theorem~\ref{T:Diff} is (globally) quasisymmetric or is a restriction of a homeomorphism $F\colon\overline\Omega\to\overline{\tilde\Omega}$ that is quasiconformal in $\Omega$, then we have the following quantitative control for the bi-Lipschitz constant. Let $p_1, p_2, p_3$ be a triple of pairwise distinct points on $\partial\Omega$ that are in the positive order, i.e., when we travel along $\dee\Omega$ starting from $p_1$ and so that $\Omega$ stays to the left, we first encounter $p_2$ and then $p_3$. Assume that for some
$\delta,\sigma>0$ we have 
$$
{\rm dist}(p_i,p_j), {\rm dist}(f(p_i),f(p_j))\geq \delta,\ i\neq j, \ {\rm and}\ {\rm diam}(\Omega), {\rm diam}(\tilde\Omega)\leq \sigma.
$$ 
Then for every $p\in S$ such that 
$$
{\rm dist}(p,\partial\Omega), {\rm dist}(f(p),\partial\tilde\Omega)\geq d>0, 
$$
there exist $r>0$ and  $L\geq 1$ that depend only on $\Omega, S, \delta, \sigma$, and $d$, so that $f$ is $L$-bi-Lipschitz in $B(p,r)\bigcap S$. In addition, there exist $r'>0$ and  $L'\geq1$ that depend only on $\Omega,\tilde\Omega, S, \tilde S, \delta, \sigma$, and $d$, so that the derivative $f'$, defined by~(\ref{E:Der}), is $L'$-Lipschitz in $B(p,r')\bigcap S$.
\end{theorem}

As Lemma~\ref{L:Qcext} below shows, a quasisymmetric map between relative Schottky sets in domains $\Omega$ and $\tilde \Omega$ can be extended to a quasiconformal map between these domains. 
See Section~\ref{S:Qc} for the definition of quasiconformality.
Note that quasiconformal maps may not be differentiable at points of a set of measure zero
and they may change the Hausdorff dimension of such a set. Nevertheless, the following is an immediate corollary to Theorem~\ref{T:Diff}. 
\begin{corollary}\label{C:Hdim}
Let $S$ be a relative Schottky set of measure zero in a Jordan domain $\Omega$. Then a locally quasisymmetric orientation preserving map $f$ from $S$ onto any other relative Schottky set preserves the Hausdorff dimension of $S$. 
\end{corollary}

It is tempting to speculate, based on Theorem~\ref{T:Diff}, that the locally quasisymmetric map $f$ in the statement must be the restriction of a conformal map between the domains $\Omega$ and $\tilde\Omega$. However this can only be possible in the case when $f$ is the restriction of a M\"obius transformation. Indeed, if $f$ were the restriction of a conformal map $g\:\Omega\to\tilde\Omega$, then $g$ would map the discs bounded by the  peripheral circles of $S$  to the discs bounded by the peripheral circles of $\tilde S$. This implies that $g$ is the restriction of a M\"obius transformation to each such disc. Since $g$ is conformal, these M\"obius transformations patch together to a global M\"obius transformation. 
Theorem~\ref{T:Blu} below shows that relative Schottky sets as in Theorem~\ref{T:Diff} are not rigid with respect to locally quasisymmetric maps.

In~\cite{mB04} M.~Bonk gives the  following quasisymmetric uniformization. Let $S$ be a set in the plane homeomorphic to the standard Sierpi\'nski carpet. If the complementary components of $S$ are bounded by uniform quasicircles and are $\delta$-relatively separated from each other (see Section~\ref{S:Loew} for the definition) for some $\delta>0$, then $S$ is quasisymmetric to a Schottky set. Here we prove the following theorem. See Section~\ref{S:Tm} for the definition of $(\delta,m)$-Loewner.
\begin{theorem}\label{T:Blu}
Let $\Omega$ and $\tilde \Omega$ be Jordan domains in $\C$. Let $S$ be a relative Schottky set in $\Omega$. 
Then there exists a relative Schottky set $\tilde S$ in $\tilde \Omega$ and an orientation preserving  homeomorphism $f\: S\to \tilde S$ that is locally bi-Lipschitz. 

More precisely, let $p_i\in\partial\Omega$ and $\tilde p_i\in\partial\tilde\Omega,\ i=1,2,3$, be two triples of distinct points  in positive order. Assume that for some $\delta>0$ we have ${\rm dist}(p_i,p_j)\geq\delta,\ i\neq j$, and
$$
\max\{{\rm dist}(\tilde p,{\rm arc}(\tilde p_i, \tilde p_j))\colon i\neq j\}\geq \delta,\quad {\rm for\ all}\ \tilde p\in\partial\tilde\Omega, 
$$
where ${\rm arc}(\tilde p_i, \tilde p_j)$ denotes the arc of $\partial\tilde\Omega$ between $\tilde p_i$ and $\tilde p_j$ that does not contain $\tilde p_k$ for $k\neq i,j$. Also assume that ${\rm diam}(\Omega),{\rm diam}(\tilde \Omega)\leq\sigma$ for some $\sigma>0$. Finally assume that $\tilde\Omega$ is $(\tilde\delta,\tilde m)$-Loewner for some $\tilde m=\tilde m(\tilde\delta)>0$. Then there exists a map $f$ as above, such that for every $p\in S$ with 
$$
{\rm dist}(p,\partial\Omega)\geq d>0, 
$$
there exist $r>0$ and $L\geq 1$ that depend only on $\Omega, S, \delta,\sigma, \tilde m$, and $d$, so that $f$ is $L$-bi-Lipschitz in $B(p,r)\bigcap S$. 
\end{theorem}

The necessity of the separation condition for $\tilde p\in\dee\tilde\Omega$ in the statement above  can be demonstrated by a domain $\tilde\Omega$ that is a slight fattening of a tripod in the plane. Lemma~\ref{L:Elemmodlem} below shows that if $\tilde\Omega$ is a fixed Jordan domain, it is automatically $(\tilde \delta, \tilde m)$-Loewner for some $\tilde m=\tilde m(\tilde \delta)>0$.

The question of rigidity with respect to local quasisymmetric maps for relative Schottky sets in the unit disc is addressed by Theorem~\ref{T:Circrig} below. An analogous result for quasisymmetric maps is proved in~\cite{BKM07}.
The proof that we give here uses a completely different method.

\begin{theorem}\label{T:Circrig}
Suppose that $S$ and $\tilde S$ are relative Schottky sets of measure zero in the unit disc $\mathbb U^2$. Let $f\: S\to \tilde S$ be a locally quasisymmetric orientation preserving homeomorphism. 
Then $f$ is the restriction to $S$ of a M\"obius transformation.  
\end{theorem}

The paper is organized as follows. Section~\ref{S:Rss} provides basic results about relative Schottky sets. Sections~\ref{S:Loew} and~\ref{S:Qc} contain basic definitions and facts about Loewner spaces and quasiconformal maps, respectively.  
Our main tool in obtaining quantitative estimates is the transboundary modulus introduced by O.~Schramm, and it is discussed in Section~\ref{S:Tm}. Uniform properness of conformal maps between the interiors of relative Schottky sets with finitely many peripheral circles is addressed in Section~\ref{S:Top}, and geometric properties of such maps are established in Section~\ref{S:Geom}, after the discussion of the fixed point index in Section~\ref{S:Fp}. Section~\ref{S:Anal} deals with analytic properties of quasisymmetric maps between relative Schottky sets. 
Theorem~\ref{T:Crit} is proved in Section~\ref{S:Ex}.
Locally bi-Lipschitz uni\-for\-mi\-za\-tion is presented in Section~\ref{S:Lblu}, and the rigidity with respect to local quasisymmetric maps in Section~\ref{S:Rig}. Section~\ref{S:Proof} contains a proof of Theorem~\ref{T:Diff} and Section~\ref{S:Quantlip} contains a proof of Theorem~\ref{T:Quantlip}.

\medskip
\noindent
{\bf Acknowledgment.} The author is grateful to Mario Bonk for numerous conversations that inspired this work and for his carefully reading the manuscript. In fact, the conclusion of differentiability in Theorem~\ref{T:Diff} is his conjecture, stated after the local bi-Lipschitz property was established using a somewhat different technique. Also, the earlier version of this paper contained a certain relative separation assumption that is unnatural and complicates applications. The author was able to remove it after discussions with Mario. The author is also grateful to the anonymous referee for many comments and suggestions that helped to  substantially improve the presentation. 

\section{Relative Schottky sets}\label{S:Rss}

\no
Recall that if $S$ is a relative Schottky set in a domain $\Omega\subseteq\C$, we write $S=\Omega\setminus\bigcup_{i\in I}B_i$, where $B_i$ are open discs with pairwise disjoint closures $\overline B_i,\ i\in I$, that are contained in $\Omega$. The family of discs $\{B_i\}_{i\in I}$, necessarily countable, is uniquely determined by $S$ as the collection of complementary components of $S$ in $\Omega$. 
If $I$ is finite, we call the interior of $S$ a \emph{relative circle domain}, following~\cite{HS95a}, \cite{HS97}.

When we speak of a relative Schottky set in a domain $\Omega$ in $\Sph^2$, we assume that $\Omega\neq \Sph^2$.
Identifying $\Sph^2$ and $\C\bigcup\{\infty\}$, 
we conclude that there is no loss of generality to assume that relative Schottky sets are contained in the plane $\C$. Moreover, the spherical and the Euclidean metrics in planar domains are conformally equivalent, and if the domains are bounded, they are bi-Lipschitz equivalent. Therefore we may assume that 
relative Schottky sets in bounded domains in $\C$ are endowed with the restriction of the Euclidean metric. We denote the distance between two points $p$ and $q$ in this metric by $|p-q|$.

By a \emph{curve} $\gamma$ in a topological space $X$ we mean a continuous image into $X$ of $[a,b], [a,b), (a,b]$, or $(a,b)$. If $\lim_{t\to a}\gamma(t)$ and $\lim_{t\to b}\gamma(t)$ exist, they are called the \emph{end points} of $\gamma$. We say that a curve $\gamma$ \emph{connects} two sets $E$ and $F$ if one of its end points is in $E$ and the other in $F$.

The following lemma will be used repeatedly.
\begin{lemma}\label{L:Curves}
Let $S$ be a relative Schottky set (or a relative circle domain) in a domain $\Omega$ in the plane and $l$ be a curve in $\Omega$ with end points $p, q\in \overline S$. Then there exists a curve $l'$ in $S$ 
whose end points are $p$ and $q$. 
In particular, $S$ is connected.
Moreover, if $l$ is rectifiable, then  there exists such a rectifiable curve $l'$ in $S$ with
$$
{\rm length}(l')\leq\pi\, {\rm length}(l),
$$ 
if $S$ is a relative Schottky set, and for every $\epsilon>0$,
$$
{\rm length}(l')\leq(\pi+\epsilon){\rm length}(l),
$$ 
if $S$ is a relative circle domain.
\end{lemma}
\smallskip\no
\emph{Proof.}
First let $S$ be a relative Schottky set.
We enumerate the peripheral circles $\{\dee B_i\}_{i\in I}$ of $S$ in the order of decreasing radii and proceed inductively as follows. If $l$ does not intersect any of the open discs $B_i$, we are done. 
Otherwise let $\dee B_{i_1}$ be the first peripheral circle in the list such that $l\bigcap B_{i_1}\neq\emptyset$. Let $p_{1}$ and $q_{1}$ be the first and the last points of $l$, respectively, that belong to $\dee B_{i_1}$, and let $l_{1}$ denote the part of $l$ between these points. The curve $l_{1}$ is the restriction of $l$ to some closed interval $[a_1,b_1]$. We replace $l_{1}$ by the shortest arc $l_{1}'$ of $\dee B_{i_1}$ with the same end points. The resulting curve obtained from $l$ by such a replacement of its part is denoted by $l_{p,q,1}$. If $l$ is rectifiable, it is clear that 
$$
{\rm length}(l_{p,q,1})={\rm length}(l)-{\rm length}(l_{1})+{\rm length}(l_{1}')
$$ 
and 
$$
{\rm length}(l_{1}')\leq\pi{\rm length}(l_{1}).
$$

Let $\dee B_{i_2}$ be the next peripheral circle in the list such that $l_{p,q,1}\bigcap B_{i_2}\neq\emptyset$, if it exists. We have necessarily $i_2>i_1$. Let $p_{2}$ and $q_{2}$ be the first and the last points of $l_{p,q,1}$, respectively, that belong to $\dee B_{i_2}$, and let $l_{2}$ denote the part of $l_{p,q,1}$ between these points. Clearly $l_{2}$ is also a subcurve of $l$, i.e., a restriction of $l$ to some interval $[a_2,b_2]$, and $[a_1,b_1]\bigcap[a_2,b_2]=\emptyset$.
We replace $l_{2}$ by the shortest arc $l_{2}'$ of $\dee B_{i_2}$ with the same end points and denote by $l_{p,q,2}$ the curve obtained from $l_{p,q,1}$ after such a replacement.  If $l$ is rectifiable, we have 
$$
{\rm length}(l_{p,q,2})={\rm length}(l_{p,q,1})-{\rm length}(l_{2})+{\rm length}(l_{2}'),
$$ 
and also
$$
{\rm length}(l_{2}')\leq\pi{\rm length}(l_{2}).
$$
Continuing inductively we obtain a curve $l_{p,q,k}$ 
that is disjoint from $B_i$ for all $i\leq i_k$, and if $l$ is rectifiable,
$$
{\rm length}(l_{p,q,k})={\rm length}(l)-\sum_{j=1}^k{\rm length}(l_{j})+\sum_{j=1}^k{\rm length}(l_{j}')
$$
with
$$
{\rm length}(l_{j}')\leq\pi{\rm length}(l_{j}),\quad j=1,2,\dots,k.
$$
Thus for a rectifiable $l$ we have 
$$
\begin{aligned}
{\rm length}(l_{p,q,k})&\leq\pi\bigg({\rm length}(l)-\sum_{j=1}^k{\rm length}(l_{j})\bigg)+\pi\sum_{j=1}^k{\rm length}(l_{j})\\
&=\pi{\rm length}(l),\quad k=1,2,\dots.
\end{aligned}
$$

Since $i_k\to\infty$ as $k\to\infty$ and the radii of the discs $B_i$ go to 0, the curves $l_{p,q,k}$ converge to a curve $l'$ in $S$. Moreover, if $l$ is rectifiable, the length of $l'$ is at most $\pi{\rm length}(l)$ since this is true for each $l_{p,q,k}$. We are done in the case $S$ is a relative Schottky set.

If $S$ is a relative circle domain, our construction of $l'$ terminates after finitely many steps and $l'$ is a curve in $\overline S\setminus\dee \Omega$. Then for every $\epsilon>0$ we can modify $l'$ in the neighborhoods of the finitely many arcs $l_{j}'$ to obtain a curve in $S$ with the desired properties. 
\qed

\medskip

The following lemma is elementary.
\begin{lemma}\label{L:Elementary}
Let $B(p,r)$ and $B(p',r')$ be two discs in the plane and assume that their boundary circles intersect in two points. Also assume that the arc of the boundary circle of $B(p',r')$ between the intersection points that has shorter length is contained in the complement of $B(p,r)$. Then $p'\in B(p,r)$ and $r'\leq r$. 
\end{lemma}

The proof of Lemma~\ref{L:Curves} and Lemma~\ref{L:Elementary} imply the following corollary.
\begin{corollary}\label{C:Circles}
Let $S$ be a relative Schottky set (or a relative circle domain) in a domain $\Omega$, and $B(p,r)$ be a disc such that $B(p,2r)$ is contained $\Omega$. Then for every two points $p,q\in \overline{S\bigcap B(p,r)}$ there exists a curve $l$ connecting them in $\overline{S}\bigcap B(p,2r)$,  such that ${\rm length}(l)\leq\pi|p-q|$.  
\end{corollary} 

\begin{proposition}\label{P:Topcharpercir}
 Let $S$ be a relative Schottky set in a domain $\Omega\subseteq\C$ and $C$ be a topological circle embedded in $S$.
Then $S\setminus C$ is connected if and only if $C$ is a peripheral circle of $S$. 
\end{proposition}
\smallskip\no
\emph{Proof.}
If $C=\dee B_i$ is a peripheral circle, then $S\setminus C$ is  connected by Lemma~\ref{L:Curves}, because $S'=S\setminus C$ is a relative Schottky set in $\Omega'=\Omega\setminus\overline B_i$. 

Now assume that $S\setminus C$ is connected. By the Jordan Curve Theorem, $\C\setminus C$ consists of two connected components $D_1$ and $D_2$, and therefore $S\setminus C$ belongs to one of them, say $D_1$. Since $C$ is embedded in $S$, the boundary $\dee \Omega$, and hence the complement of $\Omega$, must also belong to $D_1$. Thus $D_2$ consists of the union of discs bounded by peripheral circles of $S$. If this union consisted of more than one disc, then $D_2$ would contain a point in $S$, which is impossible. Thus $D_2$ coincides with a disc $B_i$ bounded by  a peripheral circle of $S$, and hence $C=\dee D_2=\dee B_i$ is a peripheral circle. 
\qed

\medskip

\begin{corollary}\label{C:Percir}
 If $f\: S\to\tilde S$ is a homeomorphism between relative Schottky sets in planar domains, then the image under $f$ of every peripheral circle of $S$ is a peripheral circle of $\tilde S$.
\end{corollary}

\section{Loewner spaces}\label{S:Loew}

\no
Let $(X,d,\mu)$ be a metric measure space, $\mu$ is Borel regular. If $\Gamma$ is a curve family in $X$ and $p\geq 1$, the $p$-\emph{modulus} of $\Gamma$ is
$$
{\rm Mod}_p(\Gamma)=\inf \int_X\rho^pd\mu,
$$
where the infimum is over all non-negative measurable functions $\rho$ defined on $X$, such that 
$$
\int_\gamma\rho\, ds\geq1,\quad {\rm for\ all}\quad \gamma\in\Gamma.
$$
If $X$ has Hausdorff dimension $n>1$, the $n$-modulus of a curve family $\Gamma$ is called the \emph{conformal modulus} of $\Gamma$, denoted ${\rm Mod}(\Gamma)$. For two sets $E$ and $F$ in $X$ we denote by ${\rm Mod}(E,F)$ the conformal modulus of the family of curves connecting $E$ and $F$.

If $X$ is a metric space and $E, F\subseteq X$ are two sets with positive diameters, we define the 
\emph{relative distance} between them to be 
$$
\Delta(E,F)=\frac{{\rm dist}(E,F)}{\min\{{\rm diam}(E),
{\rm diam}(F)\}}.
$$
We say that two sets $E$ and $F$ are $\delta$-\emph{relatively separated}, $\delta>0$, if $\Delta(E,F)\geq\delta$. 
The importance of the relative distance stems from the fact that for Ahlfors regular Loewner metric measure spaces (see below for the definitions), examples of which include $\Sph^n, \R^n$, and the unit ball $\mathbb U^n,\ n\geq2$, it gives a quantitative control for the conformal modulus of the family of curves connecting  the given sets, see~\cite{jH01}, \cite{HK98}.

A path-wise connected metric measure space $(X,d,\mu)$ of Hausdorff dimension $n>1$ is called a \emph{Loewner space} if there exists a decreasing function $\phi\: (0,\infty)\to(0,\infty)$ such that 
$$
{\rm Mod}_n(E,F)\geq \phi(t)
$$
for all $E,F\subseteq X$ disjoint continua with 
$$
\Delta(E,F)\leq t.
$$

A metric measure space $(X,d,\mu)$ is called \emph{Ahlfors $n$-regular} if there exists a constant $C\geq1$ such that
$$
\frac1Cr^n\leq\mu(\overline{B(p,r)})\leq Cr^n,
$$ 
for every $p\in X$ and $0<r\leq{\rm diam}(X)$, where $\overline{B(p,r)}$ denotes the closure of ${B(p,r)}$. A metric space $(X,d)$ is said to be \emph{linearly locally connected} if there exists a constant $C\geq1$ such that for every $p\in X$ and $r>0$, every pair of points in $B(p,r)$ can be joined by a continuum in $B(p,Cr)$, and every pair of points in $X\setminus B(p,r)$ can be joined by a continuum in  $X\setminus  B(p,r/C)$. 

If $(X,d,\mu)$ is an Ahlfors $n$-regular Loewner space, such as $\Sph^n, \R^n$, or $\mathbb U^n,\ n\geq2$, then $X$ is linearly locally connected, the function $\phi$ above can be chosen to be a homeomorphism, and there exists a decreasing homeomorphism $\psi\: (0,\infty)\to(0,\infty)$ such that
$$
{\rm Mod}_n(B(p,r),E)\leq \psi(\Delta(B(p,r), E)),
$$
for every $p\in X, r>0$, and $E\subseteq X$ a continuum disjoint from $B(p,r)$.

\section{Quasiconformal maps}\label{S:Qc}

\no
Let $F\: X\to\tilde X$ be a homeomorphism between two metric spaces $(X, d_X)$ and $(\tilde X,d_{\tilde X})$. 
The \emph{dilatation} of $F$ at $p\in X$ is defined by
\begin{equation}\label{E:Qc}
H_F(p)=\limsup_{r\to 0+}\frac{L_F(p,r)}{l_F(p,r)},
\end{equation}
where 
\begin{equation}\label{Lf}\notag
\begin{aligned}
L_F(p,r)&=\sup\{d_{\tilde X}(F(p),F(q))\co\ q\in X,\ d_{X}(p,q)\leq r\},\quad{\rm and}\\
l_F(p,r)&=\inf\{d_{\tilde X}(F(p),F(q))\co\ q\in X,\ d_{X}(p,q)\geq r\}. 
\end{aligned}
\end{equation}
The map $F$ is called \emph{quasiconformal} if 
$$
\sup_{p\in X}H_F(p)<+\infty.
$$
A quasiconformal map $F\: X\to\tilde X$ is called $H$-\emph{quasiconformal}, if 
$$
H_F(p)\leq H\quad {\rm for\ every}\ p\in X.
$$
If $H=1$, the map $F$ is called \emph{conformal}. 

It is immediate that every $\eta$-quasisymmetric map is $H$-qua\-si\-con\-for\-mal for $H=\eta(1)$. The converse holds for Ahlfors regular Loewner spaces, 
see~\cite{HK98}. Namely, suppose that $X$ and $\tilde X$ are Ahlfors $n$-regular metric measure spaces, $n>1$,  $X$ is a Loewner space, and $\tilde X$ is linearly locally connected. Let $f$ be an $H$-quasiconformal map from $X$ to $\tilde X$. If $X$ and $\tilde X$ are bounded spaces, then $f$ is $\eta$-quasisymmetric.
If $X$ and $\tilde X$ are unbounded and $f$ maps bounded sets to bounded sets, then $f$ is $\eta$-quasisymmetric. In both cases $\eta$ depends on $H$ and the data of $X$ and $\tilde X$.  

An orientation preserving homeomorphism $F$ between two domains in $\C$ is  quasiconformal if and only if  $F$ is  absolutely continuous on almost every line and there exists
$k,\ 0\leq k<1$, such that  
$$
|F_{\overline z}|\leq k |F_z|
$$
for almost every $z$, where 
$$
F_{\overline z}=\frac12\bigg(\frac{\dee F}{\dee x}+i\frac{\dee F}{\dee y}\bigg),\quad F_{z}=\frac12\bigg(\frac{\dee F}{\dee x}-i\frac{\dee F}{\dee y}\bigg),\quad z=x+iy.
$$

A \emph{Beltrami coefficient} on a measurable set $E\subseteq\C$ is a measurable complex-valued function $\mu$ defined on $E$ such that
$$
{\rm ess\,sup}\{|\mu(z)|\:\ z\in E\}<1.
$$ 
If $F$ is an orientation preserving  quasiconformal map between two domains in $\C$, the quotient $F_{\overline z}/F_z$ is a Beltrami coefficient, and it is denoted by $\mu_F$. If $F$ is orientation reversing, then we define $\mu_F=\mu_{\overline F}$. 
The Measurable Riemann Mapping Theorem states that if $\mu$ is an arbitrary Beltrami coefficient in a domain $\Omega\subseteq\C$, the Beltrami equation
$$
F_{\overline z}=\mu(z) F_{z}
$$
has an orientation preserving  quasiconformal solution $F$. 
See~\cite{lA66} and \cite{LV} for these and other facts about quasiconformal mappings.

Each Beltrami coefficient $\mu$ on a measurable set $E\subseteq\C$ defines a conformal class of measurable Riemannian metrics $ds^2$ on $E$ by
$$
ds^2=\lambda(z)|dz+\mu(z)d\overline z|^2,
$$
where $\lambda$ is a measurable function on $E$ that is positive almost everywhere.

If $\tilde \mu$ is a Beltrami coefficient on a measurable set $\tilde E\subseteq\C$ that defines  a measurable Riemannian metric $d\tilde s^2$ and $F\:\ \Omega\to\tilde\Omega$ is an orientation preserving  quasiconformal map from a domain $\Omega$ to a domain $\tilde \Omega$ that contains $\tilde E$, then there exists a well-defined  pull-back measurable Riemannian metric $ds^2=F^*(d\tilde s^2)$ on $E=F^{-1}(\tilde E)$, and it lies in a conformal class determined by some Beltrami coefficient $\nu$. We denote $\nu=F^*(\mu)$, and call it the \emph{pull-back} Beltrami coefficient.  

If $S$ is a Schottky set, the subgroup $G_S$ of the group of M\"obius transformations generated by reflections in the peripheral  circles of $S$ is a discrete group, and it is called a \emph{Schottky group associated to} $S$, see~\cite[Section~3]{BKM07}. The sets $m(S),\ m\in G_S$, form a \emph{measurable partition} of the set 
$$
S_{\infty}=\bigcup_{m\in G_S}m(S),
$$ 
i.e., for every two distinct elements $m_1$ and $m_2$ of $G_S$, the sets $m_1(S)$ and $m_2(S)$ intersect in a set of measure zero.

If $S$ is a positive measure Schottky set in the plane $\C$ and $\mu$ is an arbitrary Beltrami coefficient on $S$, then there exists a well-defined Beltrami coefficient $\mu_{\infty}$ on $S_{\infty}$, such that $\mu_\infty=\mu$ on $S$ and which is invariant under $G_S$, i.e.,
$$
m^*(\mu_{\infty})=\mu_{\infty}
$$
for all $m\in G_S$. This follows from the fact that $m(S),\ m\in G_S$, form a measurable partition of $S_\infty$. We extend $\mu_\infty$ to $\C\setminus S_{\infty}$ by zero, and let $F$ be a solution to the Beltrami equation 
$$
F_{\overline z}=\mu_\infty(z) F_{z}.
$$
The map $F$ is quasiconformal in the plane and it maps $S$ to a Schottky set $\tilde S$, see~\cite[Lemma~7.2]{BKM07}.

Since the Euclidean and the spherical metrics in $\C$ are conformally equivalent, a homeomorphism between two domains in $\C$ is quasiconformal in one of these metrics if and only if it is quasiconformal in the other. Thus the map $F$ above extends by $F(\infty)=\infty$ to a quasiconformal homeomorphism of $\Sph^2$, and since $\Sph^2$ is a Loewner space, it is a quasisymmetric map.  We collect  these facts in the following lemma.
\begin{lemma}\label{L:Schpos}
If $S$ is a positive measure Schottky set in the plane and $\mu$ is a Beltrami coefficient in $S$, then there exists an orientation preserving quasiconformal homeomorphism $F$ of the plane with $\mu_F=\mu$ on $S$, that maps $S$ to a Schottky set $\tilde S$. Moreover, the map $F$, extended by $F(\infty)=\infty$, restricts  to a quasisymmetric map of Schottky sets $S\bigcup\{\infty\}$ and $\tilde S\bigcup\{\infty\}$ in the sphere. 
\end{lemma}

Conformal maps are known to preserve the conformal modulus of a curve family.
Quasiconformal maps may change the conformal modulus.
The following lemma is elementary and we leave details to the reader.
\begin{lemma}\label{L:Changemod}
Let $\mathbb U^2$ be the unit disc in the plane and $z_1,z_2,z_3$, and $z_4$ be four distinct points in positive order on the boundary $\dee\mathbb U^2$, i.e., when we start from $z_1$ and travel along $\dee\mathbb U^2$ so that $\mathbb U^2$ stays to the left, we first encounter $z_2$, then $z_3$, and then $z_4$. Let $\Gamma$ be a family of all curves  in $\mathbb U^2$ or the punctured unit disc $\mathbb U^*$ with one end point in the arc of $\dee \mathbb U^2$ between $z_1$ and $z_2$, and the other in the arc between $z_3$ and $z_4$. Then the conformal modulus of $\Gamma$ is a positive real number and there exists a homeomorphism $F$ of $\overline{\mathbb U^2}$, quasiconformal in $\mathbb U^2$, such that the conformal modulus of 
$$
\tilde \Gamma= \{\tilde\gamma=F(\gamma)\:\ \gamma\in \Gamma\} 
$$
is different from that of $\Gamma$.
\end{lemma}

\section{Schramm's transboundary modulus}\label{S:Tm}

\no
Let $A$ be a finitely connected domain in the plane with boundary components $C_0,C_1,\dots, C_n$, and let $\Gamma$ be a family of curves in $\C$. A \emph{mass distribution} $\rho$ in $A$ is an assignment of a non-negative measurable function $z\mapsto\rho(z)$ on $A$ and non-negative numbers $\rho_i=\rho(C_i),\ i=0,1,\dots, n$, to $C_0,C_1,\dots, C_n$, respectively. We say that a mass distribution $\rho$ is \emph{admissible} for $\Gamma$ if 
$$
l_{\rho}(\gamma)=\int_{\gamma\bigcap A}\rho(z)|dz|+\sum_{i\: \gamma\bigcap C_i\neq\emptyset}\rho_i\geq1\quad {\rm for\ all}\quad \gamma\in\Gamma.
$$
Here the integral over $\gamma\bigcap A$ is defined for every $\gamma\in\Gamma$ such that each component of $\gamma\bigcap A$ is rectifiable. Otherwise we set it to be $\infty$.
The \emph{total mass} of a mass distribution $\rho$ is defined as
$$
{\rm mass}(\rho)=\int_A\rho^2(z)dxdy+\sum_{i=0}^n\rho_i^2.
$$
The \emph{transboundary modulus} of $\Gamma$ with respect to $A$ is defined as
$$
{\rm mod}_A(\Gamma)=\inf\{{\rm mass}(\rho)\:\rho\ {\rm is\ admissible\ for}\ \Gamma\},
$$
see~\cite{oS95}.
Recall that for the conformal modulus ${\rm Mod}(\Gamma)$ of a family of curves $\Gamma$ one only uses the \emph{mass function} $z\mapsto \rho(z)$.

It follows immediately from the definition that the transboundary modulus satisfies the following monotonicity property. 
If $\Gamma$ and $\Gamma'$ are two curve families such that every $\gamma\in \Gamma$ contains a subcurve $\gamma'\in\Gamma'$, then ${\rm mod}_A(\Gamma)\leq{\rm mod}_A(\Gamma')$.  
The transboundary modulus is also a conformal invariant. Namely, if $f$ is a homeomorphism of the plane, conformal in $A$, then 
$$
{\rm mod}_A(\Gamma)= {\rm mod}_{f(A)}(f(\Gamma)),
$$
where $f(\Gamma)=\{\tilde\gamma=f(\gamma)\: \gamma\in\Gamma\}$. 
The proof is immediate. Other, less elementary properties of the transboundary modulus are stated and proved below.

If $B$ is an open (or a closed) disc in the plane with radius $r$, and $t$ is an arbitrary positive number, we denote by $tB$ the open (or the closed) disc with the same center as $B$ and whose radius is $tr$.

The following lemma is well-known, see~\cite[Lemma~4.2]{bB87} and~\cite[Exercise~2.10]{jH01}. We give a proof for the sake of completeness.

\begin{lemma}\label{L:Boj}
Suppose that $\{B_1, B_2,\dots, B_n\}$ is a collection of disjoint open discs in the plane, $a_1, a_2,\dots, a_n$ are non-negative real numbers, and $\lambda\geq1$. Then there exists a constant $C\geq0$ that depends only on $\lambda$, such that
\begin{equation}\label{E:Boj}
\int\bigg(\sum_{i=1}^{n}a_i\chi_{\lambda B_i}\bigg)^2dxdy\leq C\sum_{i=1}^{n}a_i^2 \int\chi_{B_i}dxdy,
\end{equation}
where $\chi_E$ denotes the characteristic  function of a set $E$.
\end{lemma}
\smallskip\no
\emph{Proof.}
Let $\phi\in L^2=L^2(\R^2,dxdy)$. We denote the non-centered maximal function of $\phi$ by $M(\phi)$, i.e., 
$$
M(\phi)(x,y)=\sup_{B}\frac1{|B|}\int_{B}|\phi(s,t)|dsdt,
$$
where $B$ is an open disc containing $(x,y)$ and $|B|$ denotes its area.
Then
$$
\begin{aligned}
&\bigg|\int\sum_{i=1}^n a_i\chi_{\lambda B_i}\phi(s,t)dsdt\bigg|
=\bigg|\sum_{i=1}^na_i\int_{\lambda B_i}\phi(s,t)dsdt\bigg|\\
&\leq\sum_{i=1}^na_i\lambda^2\int_{B_i}M(\phi)(x,y)dxdy
=\lambda^2\int\sum_{i=1}^n a_i\chi_{B_i}M(\phi)(x,y)dxdy\\
&\leq\lambda^2\bigg|\bigg|\sum_{i=1}^na_i\chi_{B_i}\bigg|\bigg|_{L^2}\cdot||M(\phi)||_{L^2}
\leq H\lambda^2\bigg|\bigg|\sum_{i=1}^na_i\chi_{B_i}\bigg|\bigg|_{L^2}\cdot||\phi||_{L^2},
\end{aligned}
$$
where $H$ is an absolute constant. The last inequality is the maximal function inequality and it can be found in~\cite{eS70}.
This gives
$$
\bigg|\bigg|\sum_{i=1}^na_i\chi_{\lambda B_i}\bigg|\bigg|_{L^2}\leq H\lambda^2\bigg|\bigg|\sum_{i=1}^na_i\chi_{B_i}\bigg|\bigg|_{L^2}.
$$
Inequality~(\ref{E:Boj}) follows with $C=H\lambda^2$, since the discs $B_1,B_2,\dots,B_n$ are disjoint.
\qed

\medskip

The following lemma is elementary.
\begin{lemma}\label{L:Ell}
There exists a universal constant $C>0$ with the following property. Let $K$ be a planar continuum and $\{B_i\}_{i\in I}$ be a collection of  disjoint closed discs in the plane such that for each $i\in I$ we have
$$
B_i\bigcap K\neq\emptyset\quad{\rm and}\quad 2{\rm diam}(B_i)\geq{\rm diam }(K).
$$  
Then the cardinality of $I$ is at most $C$.
\end{lemma}

Using these lemmas we can prove that there is a uniform lower bound for the quotient of the transboundary modulus with respect to a relative circle domain to the conformal modulus, if the latter is small. 
\begin{proposition}\label{P:Modcomp}
Suppose that $\Omega$ is a planar domain and $A=\Omega\setminus\bigcup_{i\in I}\overline{B_i}$ is a relative circle domain in $\Omega$. Let $E$ and $F$ be two disjoint continua in $\overline \Omega$. 
Let $\Gamma$ be a family of curves in $\Omega$ that connect $E$ and $F$. Then there exists a universal constant $c>0$, such that 
$$
{\rm mod}_A(\Gamma)\geq\min\{1/(4C^2), c\,{\rm Mod}(\Gamma)\},
$$ 
where $C>0$ is the universal constant from Lemma~\ref{L:Ell}.
\end{proposition}
\smallskip\no
\emph{Proof.}
We may assume that ${\rm mod}_A(\Gamma)< 1/(4C^2)$.

Let $r_i$ denote the radius of $B_i,\ i\in I$.
To prove the inequality, we let 
$0<\epsilon\leq 1/(4C^2)-{\rm mod}_A(\Gamma)$ be arbitrary, and let 
$$
\rho=\{\rho(z),\rho_i\:\ z\in A,\ i\in I\}
$$ 
be an admissible  mass distribution for the transboundary modulus such that 
$$
{\rm mass}(\rho)\leq{\rm mod}_A(\Gamma)+\epsilon.
$$
We extend the function $z\mapsto\rho(z)$ by zero in the discs $B_i,\ i\in I$, and define a mass function on $\Omega$ by
$$
\rho_\Omega(z)=2\left(\rho(z)+\sum_{i\in I}\frac{\rho_i}{r_i}\chi_{2B_i}(z)\right).
$$
This mass function  is admissible for $\Gamma$. Indeed, let $\gamma$ be an arbitrary curve in $\Gamma$ and let $I_\gamma$ consists of all those $i\in I$ such that $\overline B_i\bigcap{\rm closure}(\gamma)\neq\emptyset$ and $2{\rm diam}(B_i)\geq{\rm diam}(\gamma)$. By Lemma~\ref{L:Ell}, the cardinality of $I_\gamma$ is at most $C$.  For any $i\in I$ we have $\rho_i\leq\sqrt{{\rm mass}(\rho)}\leq 1/(2C)$. Thus
$
\sum_{i\in I_\gamma}\rho_i\leq 1/2,
$
and hence 
$$
\int_{\gamma\bigcap A}\rho(z)|dz|+\sum_{{i\in I\setminus I_\gamma\: \gamma\bigcap \dee B_i\neq\emptyset}}\rho_i\geq \frac12.
$$
It is clear that for $i\in I\setminus I_\gamma$ such that $\dee B_i\bigcap\gamma\neq\emptyset$, the curve $\gamma$ is not contained in $2B_i$. 
Therefore
$$
\begin{aligned}
l_{\rho_\Omega}(\gamma)&=2\int_{\gamma}\rho_\Omega(z)|dz|=2\int_{\gamma}\bigg(\rho(z)+\sum_{i\in I}\frac{\rho_i}{r_i}\chi_{2B_i}(z)\bigg)|dz|\\
&\geq2\left(\int_{\gamma\bigcap A}\rho(z)|dz|+\sum_{i\in I\setminus I_\gamma\: \gamma\bigcap \dee B_i\neq\emptyset}\frac{\rho_i}{r_i}\int_{\gamma\bigcap2B_i}|dz|\right)\\
&\geq 2\left(\int_{\gamma\bigcap A}\rho(z)|dz|+\sum_{i\in I\setminus I_\gamma\: \gamma\bigcap \dee B_i\neq\emptyset}\rho_i\right)\geq1.
\end{aligned}
$$ 
It remains to estimate the total mass of $\rho_\Omega$ in terms of the total mass of $\rho$:
$$
\begin{aligned}
{\rm mass}(\rho_\Omega)&=\int_{\Omega}\rho_\Omega(z)^2dxdy\\ &\leq 8\bigg(\int_{A}\rho(z)^2dxdy+\int_{\Omega}\bigg(\sum_{i\in I}\frac{\rho_i}{r_i}\chi_{2B_i}(z)\bigg)^2dxdy\bigg)\\
&\leq8\bigg(\int_A\rho(z)^2dxdy+C'\sum_{i\in I}\frac{\rho_i^2}{r_i^2}\int_{\Omega}\chi_{B_i}(z)dxdy\bigg)\\
&\leq8\max\{1, {C'\pi}\}{\rm mass}(\rho)\\
&\leq 8\max\{1, {C'\pi}\}({\rm mod}_A(\Gamma)+\epsilon).
\end{aligned}
$$
The second inequality is an application of Lemma~\ref{L:Boj}; the constant $C'$ is universal.
Since $\epsilon$ is arbitrary, we conclude that 
$$
{\rm Mod}(\Gamma)\leq8\max\{1, {C'\pi}\}{\rm mod}_A(\Gamma).
$$
\qed

\medskip

The following results in this section will be needed to establish the uniform properness of conformal maps between relative circle domains.

\begin{lemma}\cite[Proposition~8.7]{mB04}\label{L:modest}
There exists a universal constant $N\in \mathbb N$ and a function $\psi\colon (0,\infty)\to(0,\infty)$ with $\lim_{t\to\infty}\psi(t)=0$ that satisfy the following properties. Let $\{\overline B_i\colon i\in I'\}$ be a finite collection of pairwise disjoint closed discs in $\mathbb C$. Further, let $E,F$ be arbitrary disjoint continua in $\mathbb C\setminus\bigcup_{i\in I'}B_i$ with $\Delta(E,F)\geq 12$. Then there exists $I_0\subseteq I'$ with the number of elements $|I_0|\leq N$, such that for $\Omega_0=\mathbb C\setminus \bigcup_{i\in I_0}\overline B_i$ and $A=\Omega_0\setminus\bigcup_{i\in I'\setminus I_0}\overline B_i$ we have
$$
{\rm mod}_A(\Gamma)\leq\psi(\Delta(E,F)),
$$
where $\Gamma$ is the collection of all curves in $\Omega_0$ that connect $E$ and $F$.
\end{lemma}

\begin{lemma}\label{L:Elemmodlem}
Let $\Omega$ be a  Jordan domain. Let $E,F$ be two continua in $\overline \Omega$ such that for some $\delta>0$ we have ${\rm diam}(E), {\rm diam}(F)\geq\delta$. Then there exists $m>0$ that depends only on $\Omega$ and $\delta$ such that 
$$
{\rm Mod}(\Gamma)\geq m,
$$
where $\Gamma$ is the family of all curves in $\Omega$ that connect $E$ and $F$.
\end{lemma}
\noindent
\emph{Proof.}
Let $\phi$ be a Riemann map of $\Omega$ onto the unit disc $\mathbb U^2$. 
Since $\Omega$ is a Jordan domain, $\phi$ extends to a homeomorphism between the closures, also denoted by $\phi$. 
Since $\phi^{-1}$ is uniformly continuous, there exists $\tilde\delta>0$ that depends only on $\Omega, \delta$, and the choice of $\phi$,  such that ${\rm diam}(\phi(E)), {\rm diam}(\phi(F))\geq\tilde\delta$. 

Now, since $\mathbb U^2$ is Loewner, there exists $m>0$ that depends only on $\tilde\delta$, so that ${\rm Mod}(\tilde\Gamma)\geq m$, where $\tilde \Gamma$ is the family of all curves in $\mathbb U^2$ that connect $\phi(E)$ and $\phi(F)$.  
Conformal invariance of the modulus finishes the proof.
\qed

\medskip

We also need the following more general result.
\begin{lemma}\label{L:Multiconn}
Let $\Omega$ be a Jordan domain and $A$ be a fixed relative circle domain in $\Omega$. Let $E$ and $F$ be two continua in $\overline A$ such that ${\rm diam}(E), {\rm diam}(F)\geq\delta$ for some $\delta>0$. Then there exists $m>0$ that depends only on $\Omega, A$, and $\delta$, such that
$$
{\rm Mod}(\Gamma)\geq m,
$$ 
where $\Gamma$ is the family of all curves in $A$ that connect $E$ and $F$.
\end{lemma}
\noindent
\emph{Proof.}
Let $\phi$ be a Riemann map of $\Omega$ onto $\mathbb U^2$. As in the proof of Lemma~\ref{L:Elemmodlem}, we conclude that there exists $\tilde\delta>0$ such that ${\rm diam}(\phi(E)),{\rm diam}(\phi(F))\geq\tilde\delta$. The constant $\tilde\delta$ depends only on $\Omega,\delta$, and the choice of $\phi$.

Let $\tilde A=\phi(A)$. All the boundary components of $\tilde A$ are analytic curves.
Let $\gamma$ be a smooth simple (without self-intersections)  curve in $\tilde A$ that connects $\phi(E)$ and $\phi(F)$. Such a curve can be chosen so that its length is at most a constant $C>0$ that depends only on $\tilde A$. Moreover, we may assume that $\gamma$ has a tubular neighborhood $U$ of width at least $\epsilon>0$, where $\epsilon$ depends only on $\tilde A$ and $\tilde \delta$, and such that $U$ is foliated by smooth simple curves of length at most $C$ that connect $E$ and $F$. Let $\tilde\Gamma$ denote the family of all curves in this foliation. It is easy to see that there exists a constant $m>0$ that depends only on $C$ and $\epsilon$, such that
$$
{\rm Mod}(\tilde\Gamma)\geq m.
$$ 
See also~\cite[Proposition~7.1]{mB04}.
Monotonicity of the modulus gives the same lower bound for the modulus of all curves in $\tilde A$ that connect $E$ and $F$. Now we apply the conformal invariance of the modulus to obtain the desired result.
\qed

\medskip

The last two lemmas motivate the following definition. For some $\delta, m>0$, we say that a domain $\Omega$ is $(\delta, m)$\emph{-Loewner}, if for any two continua $E$ and $F$ in $\overline \Omega$ such that ${\rm diam}(E), {\rm diam}(F)\geq\delta$, we have
$$
{\rm Mod}(\Gamma)\geq m,
$$
where $\Gamma$ is the family of all curves in $\Omega$ that connect $E$ and $F$.

The following lemma shows that given a Jordan domain $\Omega$ and a fixed relative Schottky set $S=\Omega\setminus_{i\in I}B_i$ in it, a domain obtained by removing from $\Omega$ a fixed number of closed discs from the family $\{\overline B_i\}_{i\in I}$ is $(\delta,m)$-Loewner, quantitatively. The difference with Lemma~\ref{L:Multiconn} is in the fact that the constant $m$ is independent of which discs $\overline B_i$ are removed.

\begin{lemma}\label{L:lbm}
Let $\Omega$ be a Jordan domain, let $S=\Omega\setminus_{i\in I}B_i$ be a relative Schottky set, and let $N\in \mathbb N,\ \delta>0$ be given.  Let $I_0\subseteq I$ be a finite subset whose cardinality is at most $N$. Let $\Omega_0=\Omega\setminus\bigcup_{i\in I_0} \overline B_i$ and let $E, F$ be two continua in $\overline S$ with $\min\{{\rm diam}(E),{\rm diam}(F)\}\geq\delta$.  Then there exists $m>0$ that depends only on $\Omega, S, N$, and $\delta$,  such that 
$$
{\rm Mod}(\Gamma)\geq m,
$$
where $\Gamma$ is the family of all curves in $\Omega_0$ that connect $E$ and $F$.
\end{lemma}
\noindent
\emph{Proof.}
%
%
Let $I_s\subseteq I_0$ be the subset of all indices $i$ such that ${\rm diam}(B_i)<\delta/(6N)$. Since ${\rm diam}(E)\geq\delta$, there exists $p\in E$ with ${\rm dist}(p,B_i)\geq\delta/(6N)$ for all $i\in I_s$. Indeed, let $\{p_k\}$ be a maximal $\delta/(2N)$-separated subset of $E$. The inequality ${\rm diam}(E)\geq\delta$ implies that the set $\{p_k\}$ contains strictly more than  $N$ elements. For each $i\in I_s$, we have ${\rm dist}(p_k,B_i)<\delta/(6N)$ holds for at most one $k$ because $\{p_k\}$ is $\delta/(2N)$-separated. Since there are at most $N$ elements in $I_s$, the pigeon hole principle yields the desired $p\in E$.

Let $E_p$ denote the connected component of $E\bigcap \overline{B(p,\delta/(12 N))}$ containing $p$. Clearly, 
$$
{\rm diam}(E_p)\geq\delta/(12N)\quad {\rm and}\quad{\rm dist}(E_p,B_i)\geq\delta/(12 N)
$$
for every $i\in I_s$.
Let $\Gamma_p$ be the family of all curves in $\Omega$ that connect $E_p$ to $F$. Let $m'>0$ be such that 
$$
{\rm Mod}(\Gamma_p)\geq m'. 
$$
Such $m'$ exists by Legma~\ref{L:Elemmodlem} and it depends only on $\Omega, N$,  and $\delta$.

We choose $0<\delta_1<\delta/(6N)$ so small that if $i\in I_0$ and ${\rm diam}(B_i)\leq\delta_1$ (hence $i\in I_s$), then
$$
{\rm Mod}(\Gamma_i)<m'/(2N),
$$
where $\Gamma_{i}\subseteq\Gamma_p$ is the subfamily of all curves that intersect $B_i$. The constant $\delta_1$ depends only on $\delta, m'$, and $N$. This is possible because the Euclidean plane is Ahlfors regular and Loewner. 


We consider the following three cases. 
%

\noindent
\emph{Case 1.} ${\rm diam}(B_i)<\delta_1$ for each $i\in I_0$. Then $\Gamma\supseteq\Gamma_p\setminus\bigcup_{i\in I_0}\Gamma_i$ and the subadditivity of the modulus gives
$$
{\rm Mod}(\Gamma)\geq m'/2.
$$

\noindent
\emph{Case 2.} ${\rm diam}(B_i)\geq\delta_1$ for each $i\in I_0$. There are only finitely many such configurations and the conclusion follows from Lemma~\ref{L:Multiconn}. The constant $m$ in this case depends only on $\Omega, S, N$, and $\delta$.

\noindent
\emph{Case 3.} There exists $i_1\in I_0$ such that ${\rm diam}(B_{i_1})<\delta_1$ and $i_2\in I_0$ such that ${\rm diam}(B_{i_2})\geq\delta_1$. Let $I_1\subset I_0$ consist of those indices for which  ${\rm diam}(B_{i_2})\geq\delta_1$. There are at most finitely many relative circle domains $\Omega\setminus\bigcup_{i\in I_1}\overline B_i$, and thus, by Lemma~\ref{L:Multiconn}, there exists $m_1>0$ that depends only on $\Omega, S, N$, and $\delta$, such that 
$$
{\rm Mod}(\Gamma_p\setminus\bigcup_{i\in I_1}\Gamma_i)\geq m_1.
$$
We choose $0<\delta_2<\delta_1$ so small that if ${\rm diam}(B_i)<\delta_2$, then
$$
{\rm Mod}(\Gamma_i)<m_1/(2N).
$$
The constant $\delta_2$ depends only on $\delta, m_1$, and $N$. This is again possible because the Euclidean plane is Ahlfors regular and Loewner. Now we consider three cases as above applied to $I_0\setminus I_1$, and iterate this procedure. Every time we have Case~1 or Case~2 it terminates. The procedure can run at most $N$ times and therefore we have the desired estimate for ${\rm Mod}(\Gamma)$.
%
\qed

\medskip

\begin{corollary}\label{C:lbm}
Let $\Omega$ be a Jordan domain, let $S=\Omega\setminus_{i\in I}B_i$ be a relative Schottky set, and let $N\in \mathbb N,\ \delta>0$ be given.  Let $I'\subseteq I$ be a finite subset and $I_0\subseteq I'$ be a subset whose cardinality is at most $N$. Let $\Omega_0=\Omega\setminus\bigcup_{i\in I_0} \overline B_i$ and $A=\Omega_0\setminus\bigcup_{i\in I'\setminus I_0}\overline B_i=\Omega\setminus\bigcup_{i\in I'}\overline B_i$.  Let $E, F$ be two continua in $\overline S$ with $\min\{{\rm diam}(E),{\rm diam}(F)\}\geq\delta$.  Then there exists $\epsilon>0$ that depends only on $\Omega, S, N$, and $\delta$,  such that 
$$
{\rm mod}_A(\Gamma)\geq\epsilon,
$$
where $\Gamma$ is the family of all curves in $\Omega_0$ that connect $E$ and $F$.
\end{corollary}
\no
\emph{Proof.}
This follows from Lemma~\ref{L:lbm} and Proposition~\ref{P:Modcomp}.
\qed

\medskip

%
%
%
%

\begin{lemma}\label{L:udom}
Suppose that $\Omega$ is a bounded domain in the plane 
and $S=\Omega\setminus\bigcup_{i\in I} B_i$ a relative Schottky set in $\Omega$. Let $I'\subseteq I$ be a finite subset 
and $A=\Omega\setminus\bigcup_{i\in I'}\overline B_i$ be a relative circle domain in $\Omega$. 
Let $c>0$, let $z_0\in\dee\Omega$, and let $\Gamma$ be the family of curves (we assume it is non-empty) in $\Omega$ so that each curve $\gamma\in\Gamma$ has end points in 
$\C\setminus\bigcup_{i\in I'} B_i$ and connects the complementary components of 
$$
\mathcal L=\{z\:\ c/t<|z-z_0|<c\}.
$$
Then there exist a universal number $T_0>1$ and a function $\phi\colon (T_0,\infty)\to(0,\infty)$ with $\lim_{t\to\infty}\phi(t)=0$, such that
$$
{\rm mod}_A(\Gamma)\leq \phi(t),\quad t\geq T_0.
$$
The  function $\phi$ depends only on $\Omega, S$, and $c$.
\end{lemma}
\smallskip\no
\emph{Proof.}
Let $N$ and $\psi$ be the number and the function from Lemma~\ref{L:modest}. There exists $T_0$ large enough so that for all $t\geq T_0$ we have $\Delta(E,F)\geq12$ and $\psi(\Delta(E,F))\leq1$, where $E=\{|z-z_0|\leq c/t\}\bigcap(\C\setminus\bigcup_{i\in I'} B_i)$ and $F=\{c\leq |z-z_0|\leq R\}\bigcap(\C\setminus\bigcup_{i\in I'} B_i)$ for some $R>0$ so large that $B(z_0, R)$ contains $\Omega$. By Lemma~\ref{L:modest}, 
$$
{\rm mod}_A(\Gamma')\leq 1,
$$
where $\Gamma'$ consists of all curves in $\Gamma$ that avoid the closures of at most $N$ discs in the family $\{B_i\colon i\in I'\}$. Let $\rho$ be an admissible mass distribution for $\Gamma'$. Augmenting $\rho$ by assigning weight 1 to each of the $N$ discs above, we obtain an admissible mass distribution for $\Gamma$, and thus 
$$
{\rm mod}_A(\Gamma)\leq N+1.
$$
This holds for all $t\geq T_0$, independent of $c$.

We now proceed as follows. Let $c_0=c$ and
$$
\mathcal L_0=\{z\: c_0/T_0<|z-z_0|<c_0\}.
$$
Let $c_1\leq c_0/T_0$ be a positive number such that no $\dee B_i,\ i\in I$, intersects both, $\{|z-z_0|=c_0/T_0\}$ and $\{|z-z_0|=c_1\}$. Since $\dee B_i,\ i\in I$, are disjoint from $\dee \Omega$, such a number exists. It depends only on $\Omega, S$, and $c$. We next look at the annulus
$$
\mathcal L_1=\{z\: c_1/T_0<|z-z_0|<c_1\}.
$$
Arguing inductively, we obtain a decreasing sequence of numbers $(c_k)$ and a sequence of annuli $(\mathcal L_k)$, so that no $\dee B_i,\ i\in I$, intersects any two of them. The number $c_k$ and the annulus $\mathcal L_k$ depend only on $\Omega, S$ and $c$.

Let $\epsilon>0$ be given. We choose $T$ so large that for $t\geq T$ the annulus $\mathcal L$ contains the first $n=\left[{(N+2)}/{\epsilon}\right]+1$ annuli $\mathcal L_k$. From above, we know that for each annulus
$\mathcal L_k$
we have ${\rm mod}_A(\Gamma_k)\leq N+1$, where $\Gamma_k$ is the family of curves in $\Omega$ with end points in $\C\setminus\bigcup_{i\in I'} B_i$, and that connect the complementary components of $\mathcal L_k$. Let $\rho_k$ be an admissible mass distribution for $\Gamma_k$ such that ${\rm mass}(\rho_k)\leq N+2$. We define a mass distribution $\rho$ for $\Gamma$ by setting it to be $\rho_k/n$ on $A\bigcap \mathcal L_k$, to be $\rho_k(\dee B_i)$ for $i\in I'$ such that $\dee B_i\bigcap\mathcal L_k\neq\emptyset$,  and 0 elsewhere. The mass distribution $\rho$ is well defined as follows from our construction of the annuli $\mathcal L_k$. This is clearly an admissible mass distribution for $\Gamma$ and its mass is at most $(N+2)/n<\epsilon$.
\qed

\medskip


\section{Uniform properness}\label{S:Top}

\no
Let $\Omega$ and $\tilde \Omega$ be Jordan domains in $\C$, and let $A=\Omega\setminus\bigcup_{i=1}^n\overline{B_i}$ and $\tilde A=\tilde\Omega\setminus\bigcup_{i=1}^n\overline{\tilde B_i}$ be relative circle domains.  
Suppose that $g\: A\to {\tilde A}$ is a conformal map. Such a map extends to a homeomorphism of the closures $\overline A$ and $\overline{\tilde A}$, also denoted by $g$. We will assume throughout that  $g(\dee B_i)=\dee \tilde B_i,\ i=1,2,\dots, n$. We have necessarily $g(\dee \Omega)=\dee \tilde \Omega$. By using Schwarz reflections in circles $\dee B_i$ and $\dee \tilde B_i,\ i=1,2,\dots,n$, we can extend the map $g$ conformally in a neighborhood of each $\dee B_i,\ i=1,2, \dots, n$. 

Also, since $\Omega$ and $\tilde \Omega$ are Jordan domains, we can extend $g$ to a homeomorphism of $\C$ by extending it first in each $B_j$ and then to $\C\setminus\overline\Omega$. Indeed, we can first extend $g$ in each $B_j$ radially, i.e., if $B_j=B(z_j,r_j)$, $\tilde B_j=B(\tilde z_j,\tilde r_j)$, and $g(z_j+r_je^{i\theta})=\tilde z_j+\tilde r_je^{i\tilde\theta}$, we can define 
$$
g(z)=\tilde z_j+r\tilde r_je^{i\tilde\theta}/r_j,\quad z=z_j+re^{i\theta},\ 0\leq r\leq r_j.
$$
The extended map is a continuous one-to-one map of $\overline\Omega$ onto $\overline{\tilde\Omega}$. Since the inverse map is defined by the same formula, it is a homeomorphism.
Once $g$ is extended to all of $\overline \Omega$, we can extend it homeomorphically to $\C\setminus\overline\Omega$ as follows. 
The domains $\Omega$ and $\tilde\Omega$ are Jordan domains, and by the Jordan--Sch\"onflies theorem there exist homeomorphisms $\phi$ and $\tilde\phi$ of $\C$ that take $\Omega$ and $\tilde\Omega$ to the unit disc $\mathbb U^2$, respectively.
Let $R$ denote the reflection in $\dee\mathbb U^2$, i.e., $R(z)=1/\bar z$. 
Now we can define $g$ in $\C\setminus\overline\Omega$ by the formula
$$
g=\tilde\phi^{-1}\circ R\circ\tilde\phi\circ g\circ\phi^{-1}\circ R\circ\phi.
$$
The resulting map is still denoted by $g$. In what follows it will be clear which extension of $g$ is used.

The proposition that follows establishes  the uniform properness of the map $g^{-1}$. 

\begin{proposition}\label{P:Comp} 
Let $\Omega,\tilde\Omega$ be Jordan domains and let $S=\Omega\setminus\bigcup_{i\in I}B_i$ be a relative Schottky set in $\Omega$. Let $I'\subseteq I$ be a finite subset and $A=\Omega\setminus\bigcup_{i\in I'}\overline B_i$ be a relative circle domain. Let $g\colon A\to\tilde A$ be a conformal map, where $\tilde A$ is a relative circle domain in $\tilde \Omega$.
Further, let $p_1,p_2$, and $p_3$ be a triple of distinct points on $\partial\Omega$. Assume that for some $\delta>0$ we have 
${\rm dist}(p_i,p_j)\geq\delta,\ i\neq j$, and 
$$
\max\{{\rm dist}(\tilde p, {\rm arc}(g(p_i),g(p_j)))\colon i\neq j\}\geq\delta,\quad{\rm for\ all}\ \tilde p\in\dee\tilde\Omega,
$$ 
where ${\rm arc}(g(p_i),g(p_j))$ denotes the arc of $\partial\Omega$ between $g(p_i)$ and $g(p_j)$ that does not contain $g(p_k),\ k\neq i,j$. Let $K\subseteq\Omega$ be a compact subset and $0<\epsilon\leq{\rm dist}(K,\partial\Omega)$. Then there exists $\tilde\epsilon>0$ that depends only on $\Omega, S,\delta$, and $\epsilon$, such that
$$
{\rm dist}(g(K),\partial\tilde\Omega)\geq\tilde\epsilon.
$$
\end{proposition}
\smallskip\no
\emph{Proof.} 
We assume that $K\bigcap\overline{A}$ is not empty, otherwise the conclusion is immediate.

Let $p$ be a point in $K\bigcap\overline{A}$ and
$\tilde p=g(p)$. Let $\tilde d$ denote ${\rm dist}(\tilde p,\dee\tilde\Omega)$. We need to find a lower bound for $\tilde d$. Assume that $\tilde d\leq\delta/(4\pi)$. Let $\tilde E$ be a shortest straight line segment that connects $\tilde p$ to $\dee\tilde\Omega$. So, ${\rm length}(\tilde E)=\tilde d$. Let $\tilde q\in\dee\tilde\Omega$ denote the end point of $\tilde E$. Our assumption implies that there exists $\tilde F$, an arc between some $g(p_i)$ and $g(p_j),\ i\neq j$, that does not contain $g(p_k),\ k\neq i,j$, such that
$$
{\rm dist}(\tilde q, \tilde F)\geq\delta.
$$
Elementary geometry shows that 
$$
{\rm dist}(\tilde E, \tilde F)\geq\delta/2.
$$

By Lemma~\ref{L:Curves} there exists a curve $\tilde E'\in\overline{\tilde A}$ with the same end points as $\tilde E$ and such that ${\rm length}(\tilde E')\leq\pi\tilde d$. Then ${\rm dist}(\tilde E',\tilde F)\geq\delta/2-\pi\tilde d\geq\delta/4$, and thus
\begin{equation}\label{E:1a}
\Delta(\tilde E',\tilde F)\geq\frac{\delta/4}{\pi\tilde d}=\frac{\delta}{4\pi\tilde d}.
\end{equation}
Assume further that $\tilde d\leq\delta/(48\pi)$. Then $\Delta(\tilde E',\tilde F)\geq12$.
Applying Lemma~\ref{L:modest} we conclude that there exists a universal constant $N\in\N$ and $I_0\subseteq I'$ with cardinality $|I_0|\leq N$, such that for $\tilde\Omega=\tilde\Omega_0\setminus\bigcup_{i\in I_0}\overline{\tilde B}_i$  we have
\begin{equation}\label{E:1b}
{\rm mod}_{\tilde A}(\tilde\Gamma)\leq\psi(\Delta(\tilde E',\tilde F)),
\end{equation}
where $\tilde\Gamma$ is the collection of curves in $\tilde\Omega_0$ that connect $\tilde E'$ and $\tilde F$, and $\psi\colon (0,\infty)\to(0,\infty)$ is a function with $\lim_{t\to\infty}\psi(t)=0$.

Let $\Omega_0=\Omega\setminus\bigcup_{i\in I_0}\overline B_i$.
Let $E=g^{-1}(\tilde E')$ and $F=g^{-1}(\tilde F)$. 
Then ${\rm diam}(E)\geq{\rm dist}(K,\dee\Omega)\geq\delta$ and $F={\rm arc}(p_i,p_j)\subseteq\dee\Omega$, and thus ${\rm diam}(F)\geq\delta$.

Let $\Gamma$ be the family of curves in $\Omega_0$ connecting $E$ and 
$F$.
The invariance of the transboundary modulus gives 
\begin{equation}\label{E:1c}
{\rm mod}_A(\Gamma)={\rm mod}_{\tilde A}(\tilde\Gamma). 
\end{equation}

Corollary~\ref{C:lbm} then gives $\epsilon>0$ that depends only on $\Omega, S$, and $\delta$, such that
\begin{equation}\label{E:1d}
{\rm mod}_A(\Gamma)\geq\epsilon.
\end{equation}
Combining (\ref{E:1a}), (\ref{E:1b}), (\ref{E:1c}), (\ref{E:1d}) with the fact that $\lim_{t\to\infty}\psi(t)=0$, we obtain the desired estimate for $\tilde d$.
\qed

\medskip

The uniform properness of the map $g$ is given by the following result.

\begin{proposition}\label{P:Comp1}
Let $\Omega,\tilde\Omega$ be Jordan domains and let $S=\Omega\setminus\bigcup_{i\in I}B_i$ be a relative Schottky set in $\Omega$. Let $I'\subseteq I$ be a finite subset and $A=\Omega\setminus\bigcup_{i\in I'}\overline B_i$ be a relative circle domain. Let $g\colon A\to\tilde A$ be a conformal map, where $\tilde A$ is a relative circle domain in $\tilde \Omega$. Further, let $p_1,p_2$, and $p_3$ be a triple of distinct points on $\partial\Omega$ and assume that ${\rm dist}(p_i,p_j), {\rm dist}(g(p_i),g(p_j))\geq\delta,\ i\neq j$, for some $\delta>0$. Finally, assume that 
$\tilde\Omega$ is $(\tilde\delta, \tilde m)$-Loewner for all sufficiently small $\tilde \delta>0$ and some $\tilde m=\tilde m(\tilde\delta)>0$.
Let $\tilde K\subseteq\tilde\Omega$ be a compact subset with ${\rm dist}(\tilde K,\partial\tilde\Omega)\geq\tilde\epsilon$. Then there exists $\epsilon>0$ that depends only on  $\Omega, S, \delta, \tilde m(\tilde\delta)$, and $\tilde\epsilon$, such that
$$
{\rm dist}(g^{-1}(\tilde K),\partial\Omega)\geq\epsilon.
$$
\end{proposition}
\noindent
\emph{Proof.}
Assuming $\tilde K\bigcap\overline{\tilde A}\neq\emptyset$ as we may, let $\tilde p\in\tilde K\bigcap\overline{\tilde A}$, and let $p=g^{-1}(\tilde p)$. We need to find a lower bound for $d={\rm dist}(p,\dee\Omega)$.

Let $q\in\dee\Omega$ be the point closest to $p$.
As in the proof of Proposition~\ref{P:Comp}, we can find a curve $E$ in $\overline A$ that connects $p$ and $q$, and such that ${\rm length}(E)\leq\pi d$. The point $q$ cannot be close to all three sides of the topological triangle $\Omega$ with vertices $p_1,p_2,p_3$. Thus there exist a constant $c>0$, that depends only on $\Omega$ and $\delta$, and an arc $F$ on $\dee\Omega$ between $p_i, p_j$, such that ${\rm dist}(q,F)\geq c$. 

We assume that $d<c/(2\pi)$. Then $E$ and $F$ are separated by the annulus 
$$
\mathcal L=\{x\in\C\:\pi d<|z-q|<c/2\}.
$$
Lemma~\ref{L:udom} provides a function $\phi$ with $\lim_{d\to0}\phi(d)=0$, that depends only on $\Omega, S$, and $\delta$, and such that 
\begin{equation}\label{E:EE1}
{\rm mod}_A(\Gamma)\leq\phi(d), 
\end{equation}
where $\Gamma$ is the family of curves in $\Omega$ that connect $E$ and $F$. 

Let $\tilde E=g(E), \tilde F=g(F)$, and $\tilde\Gamma$ is the family of curves in $\tilde\Omega$ that connect $\tilde E$ and $\tilde F$. 
By the invariance of the transboundary modulus we have
\begin{equation}\label{E:EE2}
{\rm mod}_{\tilde A}(\tilde \Gamma)={\rm mod}_A(\Gamma).
\end{equation}
Since $\tilde E$ is a curve that connects $\tilde K$ to $\dee\tilde\Omega$, we have ${\rm diam}(\tilde E)\geq\tilde\epsilon$. Also, our assumption gives ${\rm diam}(\tilde F)\geq \delta$. Thus the Loewner property of $\tilde\Omega$ gives 
$$
{\rm Mod}(\tilde\Gamma)\geq \tilde m(\min\{\delta, \tilde\epsilon\}).
$$
Proposition~\ref{P:Modcomp} now implies that
\begin{equation}\label{E:EE3}
{\rm mod}_{\tilde A}(\tilde\Gamma)\geq \tilde m',
\end{equation}
where $\tilde m'>0$ depends only on $\tilde m$.
Putting (\ref{E:EE1}), (\ref{E:EE2}), and (\ref{E:EE3}) together finishes the proof.
\qed

\medskip

\begin{remark}
The Loewner property of $\tilde\Omega$ in the statement of Proposition~\ref{P:Comp1} is automatically satisfied if $\tilde \Omega$ is a fixed Jordan domain. This is the assertion of Lemma~\ref{L:Elemmodlem}. 
\end{remark}

\section{Fixed point index}\label{S:Fp}

\no
Topological facts such as the Argument Principle, the Poincar\'e-Hopf Index Theorem, or the Circle Index Lemma were used for establishing rigidity properties notably by  Z.-X.~He and O.~Schramm~\cite{HS93}, \cite{HS96}, and M.~Shiffman~\cite{mS41}.

If $F\: X\to\tilde X$ is a map between two sets in $\C$, a point $p\in X$ is called a \emph{fixed point} of $F$ if $F(p)=p$.
Let $\gamma$ be an oriented Jordan curve in $\C$ and let $F\:
\gamma\to\C$ be a continuous map without fixed points. The
\emph{index} of $F$ on $\gamma$ is the winding number with respect to the
origin of the closed curve 
$$
\{F(z)-z\:z\in\gamma\}. 
$$ 
Now suppose
that $F\: \Omega\to \C$ is continuous, where $\Omega$ is a domain in $\C$, and assume
that $p\in\Omega$ is an isolated fixed point of $F$.
The \emph{index} of $F$ at $p$ is defined as the index of the
restriction of $F$ to the boundary $\partial D$ of a closed disc
$\overline D$ contained in $\Omega$ that contains $p$ in its interior and does
not contain any other fixed points of $F$. Here $\dee D$ is
positively oriented with respect to $D$, i.e., the orientation of $\dee D$ is such that when we follow it, $D$ stays to the left. Using homotopies one can show
that the index at $p$ is independent of $D$. 
It is easy to check that if $F$ is complex analytic at an isolated fixed point $p$, then the index of $F$ at $p$ is positive.

The
following version of the Poincar\'e-Hopf Index Theorem can be
found in~\cite{HS93}.
\begin{theorem}\label{T:PoincareHopf}
Let $A\subset\C$ be a bounded domain whose boundary consists of finitely many disjoint Jordan curves oriented positively with respect to $A$, i.e., when we follow the orientation of each component, $A$ stays to the left. Let $F\:
\overline A\to \C$ be a continuous map defined on the closure of
$A$. Assume that $F$ does not have any fixed points on the
boundary $\partial A$, and has only finitely many fixed points in
$A$. Then the index of the restriction of $F$ to $\partial A$, i.e., the sum of the indices of the restriction of $F$ to each component of $\dee A$, is
equal to the sum of the indices of $F$ at all its fixed points.
\end{theorem}

We say that a Jordan curve $\gamma$ in the plane $\C$
\emph{encloses} a set $X$ if $X$ is contained in the (open) Jordan
domain in $\C$ whose boundary is $\gamma$. Another result
from~\cite{HS93} that we need is the following Circle
Index Lemma. A version of this was known to K.~L.~Strebel~\cite{kS51}.
\begin{lemma}\label{L:Circleindex}
Let $\gamma$ and $\tilde\gamma$ be Jordan curves in $\C$, positively
oriented with respect to the Jordan domains that they bound. Let
$f\: \gamma\to \tilde\gamma$ be an orientation preserving
homeomorphism.
\newline
\noindent
1. If $\gamma$ encloses $\tilde\gamma$, or $\tilde\gamma$ encloses
$\gamma$, then the index of $f$ is equal to one. 
\newline\noindent
2. If
$\gamma$ and $\tilde\gamma$ intersect in at most two points, then the
index of $f$ is nonnegative.
\end{lemma}

Proposition~\ref{P:Mobcomp} below uses Theorem~\ref{T:PoincareHopf} and Lemma~\ref{L:Circleindex} to establish a relationship between a conformal map of relative circle domains and a M\"obius transformation that coincides with it at a point up to the second order. This will be used to prove the local Lipschitz property of the derivative. First we need the following elementary lemma.
\begin{lemma}\label{L:Mobtr}
Given $p\in\C$ and a conformal map $g$ in a neighborhood of $p$, there exists a unique orientation preserving M\"obius transformation $m$ that satisfies $m(p)=g(p),\ m'(p)=g'(p)$, and $m''(p)=g''(p)$.
\end{lemma}
\smallskip\no
\emph{Proof.} 
Without loss of generality we may assume $p=0$. Let $m$ be written as
$$
m(z)=\frac{az+b}{cz+d},
$$
with $ad-bc=1$. Then the constants $a,b,c$, and $d$ satisfy the following system:
$$
\begin{cases}
 \frac{b}{d}=g(0),\\
\frac1{d^2}=g'(0),\\
\frac{-2c}{d^3}=g''(0),\\
ad-bc=1.
\end{cases}
$$ 
Note that since $g$ is assumed to be conformal, $g'(0)\neq0$.
This system has two solutions that lead to the same transformation $m$:
$$
\begin{aligned}
a&=\pm{\sqrt{g'(0)}}\bigg(1-\frac{g(0)g''(0)}{2g'(0)^2}\bigg),
\quad &&b=\pm \frac{g(0)}{\sqrt{g'(0)}},\\ c&=\mp\frac{g''(0)}{2g'(0)\sqrt{g'(0)}}, 
\quad &&d=\pm\frac1{\sqrt{g'(0)}}.
\end{aligned}
$$
\qed

\medskip

\begin{proposition}\label{P:Mobcomp}{\rm (cf.~\cite[Lemma~3.4]{HS96})}  Let $\Omega $ and $\tilde \Omega$ be Jordan domains in $\C$, let $A=\Omega\setminus \bigcup_{i=1}^n\overline B_i$ and $\tilde A=\tilde \Omega\setminus\bigcup_{i=1}^n\overline{\tilde B_i}$ be relative circle domains, and let $g\:\overline A\to\overline{\tilde A}$ be a homeomorphism that is conformal in $A$. Let $p$ be a point in $A$, and let $m$ be the M\"obius transformation that satisfies $m(p)=g(p),\ m'(p)=g'(p)$, and $m''(p)=g''(p)$.  Then $g(\dee \Omega)\bigcap m(\dee \Omega)\neq \emptyset$.
\end{proposition}
\smallskip\no
\emph{Proof.}
Assume for contradiction that $g(\dee \Omega)\bigcap m(\dee \Omega)= \emptyset$. Let $\psi$ be the M\"obius transformation given by
$$
\psi(z)=p+\frac1{z-p}. 
$$
We have $\psi(p)=\infty$ and $\psi\circ\psi=\id$, and use $\psi$ to replace $p$ by $\infty$. We introduce an auxiliary map
$$
h=\psi\circ m^{-1}\circ g\circ\psi.
$$
Note that here we used $\psi=\psi^{-1}$.
Since $m$ and $g$ agree at $p$ to the second order, 
$$
m^{-1}\circ g(z)=z+o((z-p)^2),\quad {\rm as}\ z\to p.
$$
For the map $h$ this gives
$$
 h(z)-z\to0,\quad {\rm as}\ z\to\infty.
$$
The map $h$ is analytic in $\psi(A)\setminus\{\infty\}$ since the only solution to $h(z)=\infty$ in $\psi(A)$ is $z=\infty$. 
Also, we can extend the function $z\mapsto h(z)-z$ analytically to a neighborhood of $\infty$ by setting the value at $\infty$ to be 0.

Since non-constant analytic functions are open maps, there exists $\epsilon>0$ such that for every $a,\ 0<|a|<\epsilon$, the equation 
\begin{equation}\label{E:Fixedpt}
h(z)-z-a=0
\end{equation}
has a solution in $\psi(A)\setminus\{\infty\}$. We can choose such $a\neq0$, sufficiently close to 0, so that $\dee\psi(\Omega)$ and its image under $z\mapsto h(z)-a$ do not intersect, and each $\dee B_i,\ i=1,2,\dots,n$, intersects its image under $z\mapsto h(z)-a$ in at most two points.  We choose a circle $\dee B$ centered at the origin with radius so large that the interior $B$ contains a solution to equation~(\ref{E:Fixedpt}) along with all the peripheral circles of $\psi(A)$, and such that the index of $z\mapsto h(z)-a$ on $\dee B$ is zero. The latter can be achieved because $a\neq0$.

Now we compute the index of the restriction of $z\mapsto h(z)-a$ to the boundary of $\psi(A)\bigcap B$, oriented positively. 
Since $\dee \psi(\Omega)$ does not intersect its image under $z\mapsto h(z)-a$, the index on $\dee\psi(\Omega)$ is non-negative according to Lemma~\ref{L:Circleindex}. 
Therefore this index is non-positive if the orientation of $\dee \psi(\Omega)$ agrees with that of the domain $\psi(A)\bigcap B$, because $\dee \psi(\Omega)$ is an interior boundary component for that domain. Likewise, our choice of $a$ ensures that the index is non-positive on each $\dee B_i,\ i=1,2,\dots,n$, when its orientation agrees with that of $\psi(A)\bigcap B$. Since the winding number on $\dee B$ is zero, we conclude that the winding number on $\dee(\psi(A)\bigcap B)$ is non-positive. However this is impossible by Theorem~\ref{T:PoincareHopf}, because the map $z\mapsto h(z)-a$ is analytic in $\psi(A)\bigcap B$ and has at least one fixed point there, a solution to~(\ref{E:Fixedpt}). 
\qed

\medskip

\section{Geometric properties}\label{S:Geom}

\no
As in Section~\ref{S:Top}, 
let $\Omega$ and $\tilde \Omega$ be Jordan domains in $\C$
and let $A=\Omega\setminus\bigcup_{i=1}^n\overline{B_i}$ and $\tilde A=\tilde\Omega\setminus\bigcup_{i=1}^n\overline{\tilde B_i}$ be relative circle domains in $\Omega$ and $\tilde \Omega$, respectively. 
Also, let $g\: A\to{\tilde A}$ be a conformal homeomorphism in $A$ that takes $\dee \Omega$ to $\dee \tilde \Omega$, and $g(\dee B_i)=\dee \tilde B_i,\ i=1,2,\dots, n$.

The following version of the Schwarz--Pick Lemma for relative circle domains can be found in~\cite[Lemma~0.6]{HS93}. Versions of the Schwarz-Pick Lemma for circle packings are contained in~\cite{BS91}, \cite{bR87}, \cite{bR89}.
\begin{lemma}\label{L:SP}
Let $\mathbb U^2\subseteq \C$ be the open unit disc and assume that $\tilde\Omega\subseteq \mathbb U^2\subseteq \Omega$. Then $g$ is a contraction in the hyperbolic metric in the sense that if $p, q\in \overline A\bigcap \mathbb U^2$, then
$$
d_{\rm hyp}(g(p), g(q))\leq d_{\rm hyp}(p,q),
$$ 
where $d_{\rm hyp}$ denotes the distance in the hyperbolic metric of $\mathbb U^2$.  In particular, if $0\in \overline A$, then $|g'(0)|\leq 1$.
\end{lemma}

The proposition that follows proves the uniform local Lipschitz properties for $g, g^{-1}, g'$, and $(g^{-1})'$. Such properties for maps between circle packings were established by Z.-X.~He and O.~Schramm~\cite[Lemma~4.2 and Lemma~4.3]{HS96} in order to prove that the maps between circle packings converge to the conformal map between the domains as the sizes of circles go to zero. Our situation is different in that circles do not touch and do not degenerate.
\begin{proposition}\label{P:Lbilip}
Assume that for some $\delta>0$, a triple of points $p_1,p_2$, and $p_3$  on the boundary of $\Omega$ and their images under $g$ satisfy 
the assumptions as in Proposition~\ref{P:Comp}. Also assume that ${\rm diam}(\tilde\Omega)\leq\sigma$ for some $\sigma>0$.
Then for every $p_0\in\overline A\setminus \dee\Omega$ and $r\leq{\rm dist}(p_0,\partial \Omega)$, the map $g$ is $L$-Lipschitz and the map $g'$ is $L'$-Lipschitz in $\overline{B(p_0,r/4)\bigcap A}$, i.e.,
\begin{equation}\label{E:BiLip}
|g(p)-g(q)|\leq L|p-q|,
\end{equation}
\begin{equation}\label{E:Lipder}
|g'(p)-g'(q)|\leq L'|p-q|,
\end{equation}
for every $p,q\in \overline{B(p_0,r/4)\bigcap A}$. The constant $L$ depends only on $\sigma/r$. The constant $L'$ depends only on $\Omega, S, \delta, \sigma$, and $r$.

A  similar statement holds for $g^{-1}$ and $(g^{-1})'$. Assume that $p_1,p_2$, and $p_3$ and their images under $g$ satisfy the assumptions as in Proposition~\ref{P:Comp1} for $\delta>0$.
Also assume 
that ${\rm diam}(\Omega)\leq\sigma$ and that $\tilde\Omega$ is $(\tilde\delta, \tilde m)$-Loewner for all sufficiently small $\tilde \delta>0$ and some $\tilde m=\tilde m(\tilde\delta)>0$.
Then for every $\tilde p_0\in \overline{\tilde A}\setminus\partial\tilde\Omega$ and $\tilde r\leq{\rm dist}(\tilde p_0,\partial\tilde\Omega)$, the map $g^{-1}$ is $L_{-1}$-Lipschitz and its derivative $(g^{-1})'$ is $L_{-1}'$-Lipschitz in $\overline{B(\tilde p_0,\tilde r/4)\bigcap \tilde A}$.
The constant $L_{-1}$ depends only on $\sigma/\tilde r$, and the constant $L_{-1}'$ only on 
$\Omega, S, \delta, \sigma, \tilde m$, and $\tilde r$.
\end{proposition}
\smallskip\no
\emph{Proof.} 
Let $p,q\in  \overline{B(p_0,r/4)\bigcap A}$. By Corollary~\ref{C:Circles}, there exists a curve $l_{p,q}$ connecting them in $B(p_0,r/2)\bigcap\overline{A}$, such that ${\rm length}
(l_{p,q})\leq\pi|p-q|$.

To prove~(\ref{E:BiLip}), we start by obtaining a uniform bound for $|g'|$ on $B(p_0,r/2)\bigcap \overline{A}$.  If $p\in B(p_0,r/2)\bigcap \overline{A}$, then ${\rm dist}(p,\dee\Omega)\geq r/2$. Therefore
$$
\frac1{{\rm diam}(\tilde\Omega)}(\tilde\Omega-g(p))\subseteq \mathbb U^2\subseteq\frac2{r}(\Omega-p).
$$
Now, for
$$
g_1(z)= \frac1{{\rm diam}(\tilde\Omega)}\left(g\left(\frac{r}{2}z+p\right)-g(p)\right),
$$
Lemma~\ref{L:SP} gives $|g_1'(0)|\leq1$, i.e., 
$$
|g'(p)|\leq\frac{2{\rm diam}(\tilde\Omega)}{r}\leq\frac{2\sigma}{r}.
$$
Inequality~(\ref{E:BiLip}) now follows by integrating $g'$ over $l_{p,q}$.

The same argument gives the Lipschitz property of $g^{-1}$. The upper bound for $|(g^{-1})'|$ is $2\sigma/\tilde r$.

Similarly, to prove~(\ref{E:Lipder}), it is enough to establish a uniform bound for $|g''|$ on ${B(p_0,r/2)\bigcap A}$. Indeed, since $g''$ is continuous in $\overline A\setminus\dee\Omega$, the same uniform bound would hold in $B(p_0,r/2)\bigcap \overline{A}$, and integration of $g''$ over $l_{p,q}$ would finish the argument.
Let $p$ be an arbitrary point in ${B(p_0,r/2)\bigcap A}$, and let $m$ be the orientation preserving  M\"obius transformation such that $$m(p)=g(p), m'(p)=g'(p),\ {\rm and}\ m''(p)=g''(p). 
$$
It exists by Lemma~\ref{L:Mobtr}.
We have ${\rm dist}(p,\dee \Omega)\geq r/2$, and by Proposition~\ref{P:Comp} there exists $\tilde r$, depending only on $\Omega, S,\delta$, and $r$,  such that 
$$
{\rm dist}(g(p),\dee\tilde\Omega)\geq\tilde r. 
$$
Let $\eta=\min\{r/2,\tilde r\}$. According to Proposition~\ref{P:Mobcomp} there exists $q\in\dee\Omega$ with $m(q)\in\dee\tilde\Omega$. Thus for this $q$ we have
\begin{equation}\label{E:m}
|q-p|\geq\eta,\quad |m(q)-m(p)|\geq \eta.
\end{equation}
Also, by the first part of this lemma we have
\begin{equation}\label{E:m'}
|m'(p)|=|g'(p)|\leq L,
\end{equation}
where $L$ depends only on $\sigma$ and $r$.

Let $m$ be written in the form 
$$
m(z)=\frac{az+b}{cz+d},\quad ad-bc=1. 
$$
Then we obtain
$$
m'(p)=\frac1{(cp+d)^2},\quad m''(p)=\frac{-2c}{(cp+d)^3}.
$$
According to~(\ref{E:m}),
$$
\eta\leq|m(q)-m(p)|=\frac{|q-p|}{|cp+d||cq+d|}.
$$
This gives
$$
|cq+d|\leq\frac{|q-p|}{|cp+d|\eta}.
$$
Hence
$$
|c(q-p)|\leq|cq+d|+|-cp-d|\leq\frac{|q-p|}{|cp+d|\eta}+|cp+d|.
$$
Now we divide by $|q-p||cp+d|^3$ to get 
$$
\frac{|c|}{|cp+d|^3}\leq\frac1{|cp+d|^4\eta}+\frac1{|q-p||cp+d|^2}.
$$
Since 
$$
\frac1{|cp+d|^2}=|m'(p)|,
$$
we use~(\ref{E:m}) and~(\ref{E:m'}) to conclude that 
$$
|g''(p)|=|m''(p)|=\frac{2|c|}{|cp+d|^3}\leq \frac{2L^2+L}{\eta}.
$$

Running the same argument yields the desired Lipschitz property for $(g^{-1})'$. The only difference is that one has to apply Proposition~\ref{P:Comp1} instead of Proposition~\ref{P:Comp}. The Lipschitz constant then depends on $\Omega, S, \delta, \sigma, \tilde m(\tilde\delta)$, and $\tilde r$. 
\qed

\medskip

\section{Analytic properties}\label{S:Anal}

\no
The following lemma for Schottky sets appears in~\cite{BKM07}. The proof for relative Schottky sets follows the same lines. We include a proof for the sake of completeness.
\begin{lemma}\label{L:Qcext}
Let $f\: S\to \tilde S$ be an $\eta$-quasisymmetric map between two
relative Schottky sets $S$ and $\tilde S$ in domains $\Omega$ and $\tilde\Omega$ of the complex plane $\C$, respectively.
Then $f$ extends to an $H$-quasiconformal map $F\:\Omega\to \tilde\Omega$, where
$H$ depends only on $\eta$.
\end{lemma}
\smallskip
\no 
\emph{Proof.}
Let
$S=\Omega\setminus\bigcup_{i\in\N}\overline{B_i}$ and $\tilde S=\tilde\Omega\setminus\bigcup_{i\in \N}\overline{\tilde{B_i}}$. 
By Corollary~\ref{C:Percir}, $f$ sends a peripheral circle of $S$ to a peripheral circle of $\tilde S$. We assume that $f(\dee B_i)=\dee\tilde B_i,\ i\in \N$.
Using the Ahlfors-Beurling extension~\cite{lA66} and the fact that a disc in the plane is an Ahlfors 2-regular Loewner space,
we can extend each map 
$$
f|_{\partial B_i}
\co \partial B_i \ra \partial \tilde B_i,\ i\in \N, 
$$
to 
an  $\eta_1$-quasisymmetric map  of  $\overline{B_i}$ onto   
$\overline{\tilde B_i}$, where $\eta_1$ is independent of $i$.
These maps patch together to  a homeomorphism $ F\co \Omega\ra \tilde\Omega$ 
whose restriction to $S$ agrees with $f$ and whose restriction 
to each disc ${B_i}$ is 
an $\eta_1$-quasisymmetric  map onto $
{\tilde B_i}$.

It remains to show that there exists a constant $H\ge 1$ that depends only on $\eta$, such that for every
$z\in\Omega$, the inequality~(\ref{E:Qc}) is satisfied.
Below  we write $a\lesssim b$ for two quantities $a$ and $b$ 
if there exists
a constant $C$ that depends only on
the functions $\eta$ and $\eta_1$, such that $a\leq Cb$. 

If $z$ is inside one of the peripheral circles of $S$, then~(\ref{E:Qc})
follows   from the definition of $F$ with $H=\eta_1(1)$. Thus we need to consider only the case $z \in S$.

Since $S$ is 
connected by Lemma~\ref{L:Curves}, there exists $r_0>0$ such that the circles $\dee B(z,r)$
intersect $S$ for $0<r\leq r_0$. Let $r\in(0, r_0]$ and $w\in \dee B(z,r)$ be  arbitrary. 
Since $F|_S=f$ is 
$\eta$-quasisymmetric, it is enough to show that there exist
points $v', v'' \in S \bigcap \dee B(z,r)$ with
\begin{equation}\label{E:Ineq}
|F(v'')-F(z)|\lesssim|F(w)-F(z)|\lesssim|F(v')-F(z)|.
\end{equation}
If this is true, then $L_F(z,r)/l_F(z,r)$ is bounded by a quantity 
comparable to $\eta(1)$. 

The inequalities~(\ref{E:Ineq}) are trivial if $w$ itself is in $S$, because we can choose $v'=v''=w$. Thus we assume that $w$ is not 
in $S$, i.e., it lies in an open\begin{large}                                \end{large} disc $B_i$ bounded by one of the peripheral circles $\dee B_i$ of $S$. 
Let $v'$ denotes one of the points in $\dee B (z,r)\bigcap\partial B_i$, and 
let $u'$ be the point of intersection of $\dee B_i$ and
the line  segment $[z,w]$. Since
$$|w-u'|\leq|v'-u'|,\
|u'-z|\leq|v'-z|,\ |v'-u'|\leq2|v'-z|, 
$$
the triple
$\{z, v', u'\}$ is in $S$, and the triple $\{w, v', u'\}$ is in 
$\overline B_i$, we have
\begin{align}
|F(w)-F(z)|&\leq|F(w)-F(u')|+|F(u')-F(z)| \notag\\ 
&\lesssim|F(v')-F(u')|+|F(v')-F(z)|\lesssim|F(v')-F(z)|.\notag
\end{align}
This shows the right-hand side of~(\ref{E:Ineq}). To prove the left-hand
side inequality, we choose $v''$ in the same way as $v'$, namely to be a point in the intersection $\dee B(z,r)\bigcap\dee B_i$.
We choose $u''$ to be the preimage under $F$ of the point of
intersection of the line segment $[F(z),F(w)]$ and $F(\dee B_i)$.
Again, the triple
$\{z, v'', u''\}$ is in $S$, and the triple $\{w, v'', u''\}$ 
is in  $\overline B_i$.
We need to consider two cases.
If $|u''-z|\geq\frac12 r$, then we have 
$|v''-z|\leq2|u''-z|$, and therefore
$$
|F(v'')-F(z)|\lesssim|F(u'')-F(z)|\leq|F(w)-F(z)|.
$$
If, on the other hand, $|u''-z|<\frac12 r$, then we have 
$|v''-u''|\leq3|w-u''|$, and thus
\begin{align}
|F(v'')-F(z)|&\leq|F(v'')-F(u'')|+|F(u'')-F(z)|\notag\\
&\lesssim|F(w)-F(u'')|
+|F(u'')-F(z)|=|F(w)-F(z)|.\notag
\end{align}
This completes the proof of~(\ref{E:Ineq}), and thus of~(\ref{E:Qc})
and the lemma. \qed

\medskip

It is known that quasiconformal maps send sets of measure zero to sets of measure zero. It turns out that linear combinations of quasiconformal maps also possess this property. 
\begin{lemma}\label{L:Setsofmeasurenaught}
A linear combination of quasiconformal maps defined on a domain $\Omega\subseteq\C$ sends sets of measure zero to sets of measure zero.
\end{lemma}
\smallskip\no
\emph{Proof.} 
Let $E$ be a subset of $\Omega$ of measure zero. Without loss of generality we may assume that the closure $\overline E$ is compact and is contained in $\Omega$. Let $\epsilon>0$, and let $\{B(z_i,r_i)\}_{i\in I}$ be a cover of $E$ with 
$$10r_i\leq {\rm dist}(z_i,\dee\Omega)\quad {\rm and}\quad \sum_{i\in I}r_i^2<\epsilon.
$$ 
Assuming that each disc $B(z_i,r_i)$ intersects $E$, the union $\bigcup_{i\in I}B(z_i,r_i)$ is contained in some compact set in $\Omega$. 
By a basic covering theorem, see e.g.,~\cite[p.~2]{jH01}, there exists a disjoint subfamily $\{B(z_i,r_i)\}_{i\in I_0}$ such that
$$
\bigcup_{i\in I}B(z_i,r_i)\subseteq \bigcup_{i\in I_0}B(z_i, 5r_i)\subseteq\Omega.
$$  
Now suppose that $F$ is a $H$-quasiconformal map defined on $\Omega$.
By~\cite[Theorem~11.14]{jH01}, $F$ is $\eta$-quasisymmetric in $B(z_i, 5r_i)$ with $\eta$ depending only on $H$. This, combined with~\cite[Proposition~10.8]{jH01}, gives
\begin{equation}\label{E:diamdiam}
{\rm diam}^2(F(B(z_i,5r_i)))\leq C_1{\rm diam}^2(F(B(z_i,r_i))),
\end{equation}
for $i\in I_0$, where $C_1$ depends only on $H$.
Also, since $F|_{B(z_i,5r_i)}$ is $\eta$-quasisymmetric,
\begin{equation}\label{E:diamarea}
\begin{aligned}
{\rm diam}^2(F(B(z_i,r_i)))&\leq C_2{\area}(F(B(z_i,r_i)))\\ &=C_2\int_{B(z_i,r_i)}J_F(x,y)dxdy,
\end{aligned}
\end{equation}
$i\in I_0$, where the constant $C_2$ depends only on $H$, and $J_F$ is the Jacobian of $F$. 
Combining~(\ref{E:diamdiam}) and~(\ref{E:diamarea}), for $i\in I_0$ we obtain
$$
{\rm diam}^2(F(B(z_i,5r_i)))\leq C_3\int_{B(z_i,r_i)}J_F(x,y)dxdy,
$$
where $C_3=C_1C_2$ depends only on $H$.
The set $F(E)$ is covered by the sets
in the family $\{F(B(z_i,5r_i))\}_{i\in I_0}$. Therefore its measure is not greater than
\begin{equation}\notag
\begin{aligned}
\sum_{i\in I_0}{\rm diam}^2(F(B(z_i,5r_i)))&\leq C_3 \sum_{i\in I_0}
\int_{B(z_i,r_i)}J_F(x,y)dxdy\\ &=C_3\int_{\bigcup_{i\in I_0}B(z_i,r_i)}J_F(x,y)dxdy.
\end{aligned}
\end{equation}
Since the Jacobian of a quasiconformal map is locally integrable, the last integral can be made arbitrarily small by choosing an appropriate $\epsilon$. 

Now if $F_1, F_2,\dots, F_k$ are quasiconformal maps and $\alpha_1,\alpha_2,\dots, \alpha_k$ are constants, then
\begin{equation}\notag
\begin{aligned}
&\sum_{i\in I_0}{\rm diam}^2((\alpha_1 F_1+\dots+\alpha_k F_k)(B(z_i,5r_i)))\\
&\leq C_4\sum_{i\in I_0}\sum_{j=1}^k|\alpha_j|^2{\rm diam}^2(F_j(B(z_j,5r_j)))\\
&\leq C_5\sum_{j=1}^k|\alpha_j|^2\int_{\bigcup_{i\in I_0}B(z_i,r_i)}J_{F_j}(x,y)dxdy.
\end{aligned}
\end{equation}
Here the constant $C_4$ depends only on $k$, and $C_5=C_3C_4$. Each of the integrals in the last sum can be made arbitrarily small, and thus the lemma follows.  
\qed

\medskip

\section{Proof of Theorem~\ref{T:Crit}}\label{S:Ex}

\no
From the definition of a relative Schottky set $S$ in a domain $\Omega$ in $\Sph^2$ it follows that $S'=S\bigcup(\Sph^2\setminus\Omega)$ is a Schottky set in 
$\Sph^2$ whose peripheral circles are those of $S$.

Assume that $S'$ has measure zero and let $f\: S\to \tilde S$ be a quasisymmetric map from $S$ to a relative Schottky set $\tilde S$ in a domain $\tilde \Omega$. Since $S'$ has measure zero, $S$ is dense in $S'$, and quasisymmetric maps take Cauchy sequences to Cauchy sequences, we can extend the map $f$ to a quasisymmetric map $f'$ defined on $S'$. The image $\tilde S'$ of $S'$ under $f'$ is a Schottky set in $\Sph^2$. Indeed, extending $f'$ homeomorphically in discs bounded by the peripheral circles of $S'$ we obtain a homeomorphism of $\Sph^2$ onto a subset $\tilde\Omega'$ of $\Sph^2$. An application of the Borsuk--Ulam Theorem, see, e.g., \cite[Chapter~V, Corollary 9.4]{wM91}, shows that $\tilde\Omega'$  must be all of $\Sph^2$. Thus $\tilde S'$ is a Schottky set in $\Sph^2$ whose peripheral circles are those of $\tilde S$, and $f'$ is a quasisymmetric map from $S'$ to $\tilde S'$. By Theorem~B, the map $f'$, and hence $f$, is the restriction of a M\"obius transformation. Thus $S$ is rigid with respect to quasisymmetric maps.

Now suppose that $S'$ has positive measure. We consider two cases, depending on whether $S$ is dense in $S'$ or not. 

First we assume that  $S$ is dense in $S'$. 
By Theorem~B there exists a quasisymmetric map $f'$ from $S'$ to a Schottky set $\tilde S'$ in $\Sph^2$ that is not the restriction of a M\"obius transformation. The map $f'$ extends to a homeomorphism of $\Sph^2$ by extending it in discs bounded by the peripheral circles of $S'$. The restriction  $f=f'|_S$ is a quasisymmetric map of $S$ onto a relative Schottky set, and it cannot be the restriction of a M\"obius transformation because $S$ is dense in $S'$.

Now we assume that $S$ is not dense in $S'$. We identify $\Sph^2$ with $\C\bigcup\{\infty\}$, and
without loss of generality we assume that $\Omega\subset\C$. 
Then $\C\setminus\overline\Omega$ is a non-empty open set. Let $D$ be a connected component of this set. Then $D$ is  a either a simply connected domain 
or an annulus with one boundary component at $\infty$. 
By the Riemann Mapping Theorem, see, e.g., \cite[Chapter~II, \S 2]{gG66}, there exists a  conformal map $G$ from $D$ onto the unit disc $\mathbb U^2$ 
or the punctured unit disc $\mathbb U^*$ 
in the plane. 

Let  $z_1,z_2,z_3$, and $z_4$ be distinct points in $\dee \mathbb U^2$ in positive order. 
We consider the family $\tilde \Gamma$ of curves in $\mathbb U^*$ that connect two disjoint arcs of $\dee \mathbb U^2$, one with end points $z_1$ and $z_2$, and the other with end points $z_3$ and $z_4$.    
Let $\Gamma$ be the family of curves in $D$ given by 
$$
\Gamma=\{\gamma=G^{-1}(\tilde\gamma)\:\ \tilde\gamma\in\tilde\Gamma\}.
$$
Let $H$ be a quasiconformal map defined in $\mathbb U^2$ that changes the 
conformal modulus of $\tilde \Gamma$. Such a map exists by Lemma~\ref{L:Changemod}.
By conformal invariance, ${\rm Mod}(\Gamma)={\rm Mod}(\tilde\Gamma)$, and thus $H\circ G$ changes the conformal modulus of $\Gamma$.
Let $\mu_{H\circ G}$ be the Beltrami coefficient of the quasiconformal map $H\circ G$ in $D$, and we assume that $\mu_{H\circ G}$ is extended by zero to $S'\setminus D$. 

By Lemma~\ref{L:Schpos}, there exists a quasiconformal homeomorphism $F$ of the plane such that $\mu_F=\mu_{H\circ G}$  on $S'$, and that maps $S'$ onto a Schottky set $\tilde S'$ in $\C$. The map $F$ restricts to a quasisymmetric map $f$ of $S$ to a relative Schottky set $\tilde S$ in a domain $\tilde\Omega$. 
The map $f$ cannot be the restriction of a M\"obius transformation.
If it were, then $F|_{\dee D}$ would coincide with a M\"obius transformation, and in particular $F$ would preserve the conformal modulus of $\Gamma$. But since $\mu_F=\mu_{H\circ G}$ in $D$, the map 
$$
H\circ G\circ F^{-1}\:\ F(D)\to H(\mathbb U^2)
$$
is conformal. This leads to a contradiction because 
$$
H\circ G=(H\circ G\circ F^{-1})\circ F 
$$
changes the conformal modulus of $\Gamma$.
\qed

\medskip

\section{Proof of Theorem~\ref{T:Blu}}\label{S:Lblu}

\no
The following result is contained in~\cite[Theorem~4.2]{oS96}, see also~\cite{HS95}. Previous uniformization results of this type can be found in~\cite{mB80}, \cite{aH82}.
\begin{theorem}\label{T:BHUnif}
Let $\Omega$ and $\tilde\Omega$ be Jordan domains in $\C$, and 
$A$ be a relative circle domain in $\Omega$. Let $p_i\in\dee\Omega$ and $\tilde p_i\in\dee\tilde\Omega,\ i=1,2,3$, be two triples of distinct points  in positive order. Then there exists a conformal map $g$ from $A$ to a relative circle domain $\tilde A$ in $\tilde\Omega$, whose continuous extension to $\dee A$ maps $p_i$ to $\tilde p_i,\ i=1,2,3$.    
\end{theorem}

Let $X$ be a metric space and $E, F\sub X$ be two subsets.
 The 
{\em Hausdorff distance} ${\rm dist}_H(E,F)$ is defined as the infimum of all $\epsilon
\in (0,\infty]$ such that 
$$ E \sub N_\epsilon (F) \text{ and } F \sub N_\epsilon(E), $$
where $N_\epsilon(K)$ denotes the open $\epsilon$-neighborhood of a set $K\sub X$. 
The definition immediately gives that ${\rm dist}_H(E,F)=0$ if and only if 
 $\overline E=\overline F$.

 We say that a  sequence $(A_n)$   of  sets in $X$ \emph{Hausdorff converges} to a set $S\sub X$, if 
$$ {\rm dist}_H(A_n,S)\to 0 \text{ as } n\to \infty. $$ 
The Hausdorff convergence of sets can be checked using the following simple observations. If $(A_n)$ converges to $S$,  then for each $p\in S$ there exists a sequence $(p_n)$ such that $p_n\in A_n$ and $p_n\to p$. Conversely, if  for some $p\in X$ there exist a subsequence
$(A_{n_k})$ of $(A_n)$ and corresponding points $p_{n_k}\in A_{n_k}$ with $p_{n_k} \to p,\ k\to \infty$, then $p\in \overline S$. 
In particular, this implies  that if $p\in X\setminus \overline S$, then 
$p\in X\setminus A_n$ for large $n$.

\medskip\no
\emph{Proof of Theorem~\ref{T:Blu}.}
Let $p_i\in\dee\Omega$ and $\tilde p_i\in\dee\tilde\Omega,\ i=1,2,3$, be two triples of distinct points  in positive order.
Let $\{\dee B_i\}_{i=1}^{\infty}$ be the collection of peripheral circles of $S$. 
For each $n=1,2,\dots$, we consider a relative circle domain 
$A_n=\Omega\setminus\bigcup_{i=1}^{n}\overline B_i$. Let $g_n$ be a conformal map from $A_n$ to a relative circle domain $\tilde A_n$ in $\tilde\Omega$ whose continuous extension to $\dee A$, still denoted by $g_n$, satisfies $g_n(p_i)=\tilde p_i,\ i=1,2,3$. Such a map $g_n$ is guaranteed by Theorem~\ref{T:BHUnif}.

By Propositions~\ref{P:Comp}, \ref{P:Comp1}, the maps $g_n, g_{n}^{-1},\ n=1,2,\dots$, are uniformly proper, and by Proposition~\ref{P:Lbilip}, they are uniformly locally Lipschitz. Thus a subsequence $(\tilde A_{n_j})$ Hausdorff converges to a relative Schottky set $\tilde S$ in $\tilde\Omega$. Indeed, for each peripheral circle $\dee B_{i}$ of $S$, its image $\dee \tilde B_{n,i}$ under $g_{n},\ n\geq n_i$, is a circle contained in a compact subset of $\tilde\Omega$ independent of $n$, and whose radius is uniformly bounded above and below.  
By possibly passing to subsequences and using the diagonalization argument, we may assume that for each $i$, the sequence of discs $(\tilde B_{n,i})_{n\geq n_i}$, so that $\tilde B_{n,i}$ is bounded by $\dee \tilde B_{n,i}$, converges. It is clear that the interior of  the limit is necessarily a disc, denoted by $\tilde  B_i$, and for $i\neq j,\ \tilde B_i\bigcap\tilde B_j=\emptyset$. Now let $\tilde p\in \tilde\Omega\setminus \tilde S$. There exists a subsequence $(\tilde A_{n_k})$ and, for each $k,\ \dee \tilde B_{n_k}$, a peripheral circle of $\tilde A_{n_k}$, such that $\tilde p\in B_{n_k}$, the disc bounded by $\dee \tilde B_{n_k}$. The preimage $g^{-1}_{n_k}(\dee \tilde B_{n_k})$ is eventually, i.e., for $k$ large enough, the same peripheral circle $\dee B_i$ of $S$. This easily follows from the uniform properness and the uniform bi-Lipschitz property of $g_n$. Therefore $\tilde p\in \tilde B_i$, i.e., the complement of $\tilde S$ in $\tilde\Omega$ consists of pairwise disjoint open discs. This means that $\tilde S$ is a relative Schottky set in $\tilde\Omega$. The closures of the complementary discs are pairwise disjoint as well because $g^{-1}_n$ are locally uniformly Lipschitz.

Using the Arzel\`a-Ascoli Theorem~\cite[Theorem~7.23]{wR64},  we conclude that
there exists a subsequence of $(g_{n_{j}})$ that converges locally uniformly to a continuous map $f\colon S\to\tilde S$. We can find a further subsequence $(g_{n_{j_l}})$ of $(g_{n_j})$ so that $(g^{-1}_{n_{j_l}})$ converges locally uniformly to a continuous map $h$ from $\tilde S$ to $S$.
Since $g_{n_{j_l}}\circ (g_{n_{j_l}})^{-1}=\id$ and $(g_{n_{j_l}})^{-1}\circ g_{n_{j_l}}=\id$ for all $l$, we conclude that $f$ and $h$ are homeomorphisms with $f^{-1}=h$. The map $f$ is  locally bi-Lipschitz as a limit of uniformly locally bi-Lipschitz maps.

To prove the quantitative part of this theorem, we let $p\in S$ and $d\leq{\rm dist}(p,\dee\Omega)$. Then $p\in A_n$ for every $n$.
By Proposition~\ref{P:Lbilip}, each $g_n$ is $L_1$-Lipschitz in $\overline{B(p,d/4)\bigcap A_n}$, where $L_1$ depends only on $\sigma$ and $d$. Also, by Proposition~\ref{P:Comp}, ${\rm dist}(g_n(p),\dee\tilde\Omega)\geq\tilde\epsilon$, where $\tilde\epsilon>0$ depends only on $\Omega, S, \delta$, and $d$.  
Applying Proposition~\ref{P:Lbilip} once again, we conclude that $g^{-1}_n$ is $L_2$-Lipschitz in $\overline{B(g_n(p),\tilde\epsilon/4)\bigcap\tilde A_n}$, where $L_2$ depends only on $\sigma$ and $\tilde\epsilon$. We can now choose $r$ to be $\min\{d/4,\tilde\epsilon/(4L_1)\}$ and $L=\max\{L_1,L_2\}$.
The bi-Lipschitz constant $L$ persists under taking the limits as above.
\qed

\medskip

\section{Proof of Theorem~\ref{T:Circrig}.}\label{S:Rig}

\no
By pre- and post-composing $f$ with M\"obius transformations that preserve $\mathbb U^2$ we may assume that one of the peripheral circles of $S$, say $\dee B_1$, and its image $\dee \tilde B_1$ are centered at the origin. Let $\dee B$ be an arbitrary peripheral circle of $S$, other than $\dee B_1$.
Further post-composing $f$ with a rotation and a dilation (or a \index{\footnote{\left[ }}contraction) with respect to the origin, we get a locally quasisymmetric homeomorphism $f_B$ from $S$ to a relative Schottky set in a disc $U_B$ centered at the origin, such that $f_B(\dee B)$ has the same Euclidean center as $\dee B$.
We consider two cases depending on whether all of the equalities 
\begin{equation}\label{E:Rot}
\dee U_B=\dee \mathbb U^2,\quad f_B(\dee B_1)=\dee B_1,\quad {\rm and}\quad  f_B(\dee B)=\dee B
\end{equation}
hold or not.

Assume first that there exists $B$, say $B=B_2$, such that at least one of the equalities in~(\ref{E:Rot}) fails. We will show that this case is actually impossible.
Indeed, let $D$ denote one of the following discs: the unit disc $\mathbb U^2$, or the disc  $B_1$, or the disc $B_2$. Let $\tilde D$ denote $U_{B_2}$ if $D=\mathbb U^2$, or the open disc in $\C$ bounded by $f_{B_2}(\dee D)$ in the other two cases.
There exists a M\"obius transformation $m$ such that for each choice of $D$ as above either $m(\dee\tilde D)$ is contained in  $D$ or else $m(\tilde D)$ contains $\dee D$.
This can be seen as follows. If none of the equalities in~(\ref{E:Rot}) holds, then $m=\id$ works because the corresponding circles in~(\ref{E:Rot}) have the same centers. If only one of the equalities holds, then we apply a dilation with the center at the corresponding circle and a coefficient close to one. If two of the equalities hold, then we apply a M\"obius transformation that has a repelling fixed point at the center of one of these circles and an attracting fixed point at the center of the other, with coefficients close to one. In the case when these circles are $\dee \mathbb U^2$ and $\dee B_1$, it is simply a dilation with a coefficient close to one.

Our choice of $m$ implies, in particular, that there exist constants $c>0$ and $r_0, \ 0\leq r_0<1$, such that 
\begin{equation}\label{E:Sep}
|m(f_{B_2}(z))-z|\geq c\quad {\rm  for\ all}\quad  z\in S\bigcap \{z\:\ r_0\leq|z|< 1\}.
\end{equation}
Let $r,\ r_0\leq r<1$, be chosen so that the following hold: the discs $B_1, B_2$ are contained in $\{z\:|z|\leq r\}$, there is no peripheral circle of $S$ that intersects both $\{z\: |z|=r_0\}$ and $\{z\: |z|=r\}$, and there is no peripheral circle of $S$ that has only one point of intersection with $\{z\: |z|=r\}$.
Such $r$ exists because peripheral circles of a relative Schottky set do not touch the boundary of the corresponding domain and there are only countably many of them. 
Let $C_r$ be a curve obtained from  $\{z\: |z|=r\}$ by replacing each arc inside a peripheral circle of $S$ by the arc of this peripheral circle contained in  $\{z\: |z|\leq r\}$ and with the same end points. It is a Jordan curve since peripheral circles are disjoint and their diameters go to 0. We may also assume that $r$ is chosen so that the curve $C_r$ obtained in this way does not intersect its image under $m\circ f_{B_2}$. This is still possible because $\dee \mathbb U^2\bigcap m(\dee U_{B_2})=\emptyset$ and no peripheral circle of $S$ touches $\dee \mathbb U^2$.

Let $\Omega_r$ denote the domain in $\C$ bounded by the curve $C_r$. Then $S_r=S\bigcap\Omega_r$ is a relative Schottky set in $\Omega_r$, which follows from the choice of $C_r$. The map $m\circ f_{B_2}$ is quasisymmetric in $S_r$ because $f$ is locally quasisymmetric and the closure of $\Omega_r$ is a compact subset of $\mathbb U^2$. 
Using Lemma~\ref{L:Qcext}, we can extend $m\circ f_{B_2}$ to a quasiconformal map $F$ defined in $\Omega_r$. The map $F$ also extends to a homeomorphism on the closure $\overline \Omega_r=\Omega_r\bigcup C_r$  by $F=m\circ f_{B_2}$ on $C_r$.
By Lemma~\ref{L:Setsofmeasurenaught}, the map $z\mapsto F(z)-z$ sends $S_r$ to a set of measure zero, and therefore in any neighborhood of the origin there exists a full measure subset of elements $a$ such that $F(z)-a\neq z$ for all $z\in S_r$.
In addition, an element  $a$ with this property can be chosen so close to 0 to satisfy the following.
First we may require the inequality $F(z)-a\neq z$ to hold for all $z\in C_r$. This is true for any $a$ such that $|a|<{\rm dist}(C_r,m(f_{B_2}(C_r)))$. Next, if $D$ denotes one of the domains $\Omega_r$, or  $B_1$, or $B_2$, then we may require that either $F(\dee D)-a$ is contained in  $D$, or $F(D)-a$ contains $\dee D$. This property is true for $a=0$ because of the choices of $m$ and $r$ above.  Since $D, F(D)$ are open and $\dee D, F(\dee D)$ are compact, the inclusions continue to hold for $a$ in a neighborhood of $0$.   
Finally, since the number of peripheral circles of $S_r$ is countable, an element $a$ can be chosen so that for each peripheral circle $\dee B_i$ of $S_r$, $F(\dee B_i)-a$ intersects $\dee B_i$ in at most two points.

Since 
$F_a=F-a$ does not have fixed points in $\overline S_r$, 
there are only finitely many peripheral circles of $S_r$ that enclose fixed points of $F_a$.
Now we consider a finitely connected domain $A$ obtained from $\Omega_r$ by removing $\overline B_1,\overline B_2$, and finitely many other closed discs bounded by peripheral circles of $S_r$ that enclose fixed points of $F_a$. Since $F_a$ does not have any fixed points in $A$, by Theorem~\ref{T:PoincareHopf}, the index of the restriction $F_a|_{C_r}$ is equal to the sum of the indices of the restrictions of $F_a$ to the peripheral circles of $S_r$ that are boundary components of $A$. Here $C_r$ is oriented positively with respect to $\Omega_r$ and the peripheral circles are oriented positively with respect to the discs in $\C$ that they bound. However, according to Lemma~\ref{L:Circleindex}, the indices of the restrictions of $F_a$ to $C_r, \dee B_1$, and $\dee B_2$ are equal to one, and the indices of the restrictions of $F_a$ to other peripheral circles are non-negative. This gives a contradiction.

Assume now that for every $B$ all the equalities in~(\ref{E:Rot}) hold. 
In this case $f$ must be a rotation, which can be seen as follows. 
The first equality in~(\ref{E:Rot}) implies that for every peripheral circle $\dee B$ there exists a rotation $R_B$ such that $f_B=R_B\circ f$. The middle and last equalities in~(\ref{E:Rot}) then tell us that $f(\dee B_1)=\dee B_1$ and $f(\dee B)=R_B^{-1}(\dee B)$.
Since $S$ has measure zero, every point in $S$ is an accumulation point of peripheral circles of $S$. 
Let $\dee B$ be a fixed peripheral circle of $S$, let $p$ be an arbitrary point in $\dee B$, and let $(\dee B_k)$ be a sequence of distinct peripheral circles of $S$ that accumulate at $p$. In particular, ${\rm diam}(\dee B_k)\to0$ as $k\to\infty$.   Let $R_k$ be a rotation such that $f(\dee B_k)=R_k(\dee B_k)$. Since $f$ is continuous, the sequence of rotations $(R_k)$ converges uniformly on compacta to a rotation $R$. Thus $f(p)=R(p)$, i.e., the map $f$ takes $p$ to a point that has the same distance from the origin. Assume now that $p$ is neither the closest nor the farthest point of $\dee B$ in relation to 0. Then for $f(p)$ we have two choices, the intersection points of $f(\dee B)$ with $\{z\in\C\co |z|=|p|\}$. Combined with the assumptions that  $f$ is continuous and orientation preserving, this implies that the restriction of $f$ to every peripheral circle of $S$ coincides pointwise with a rotation. 

This further implies that $f$ preserves the distances between any two points in $S$. Indeed, let $p$ and $q$ be a pair of points in $S$ and $l_{p,q}$ denote the line segment connecting them. By choosing $r$ sufficiently close to one, we can find a domain $\Omega_r$ as above that contains $l_{p,q}$. The map $f$ is quasisymmetric in $S_r=S\bigcap \Omega_r$, and therefore, by Lemma~\ref{L:Qcext}, it has a  quasiconformal extension $F_r$ to $\Omega_r$. 
Thus there is a pair of points $p'$ and $q'$ in $\Omega_r$, such that $p'$ is close to $p$, $q'$ is close to $q$, the line segment $l_{p',q'}$ connecting them is in $\Omega_r$, and $F_r$ is absolutely continuous on $l_{p',q'}$. Moreover, since $S_r$ has measure zero, using Fubini's theorem we may assume that $l_{p',q'}\subseteq \Omega_r$ spends zero length in $S_r$. 
Since the restriction of $F$ to every peripheral circle of $S$ equals $f$ and thus coincides with a rotation, $F$ maps $l_{p',q'}$ to a curve that has the same length. Therefore $f$ does not increase the distance between points. Applying the same result to $f^{-1}$, we conclude that $f$ is an isometry. 
Since $f$ preserves $\dee B_1$, it is a rotation.
\qed

\medskip

\section{Proof of Theorem~\ref{T:Diff}.}\label{S:Proof}

\no
Let $p$ be an arbitrary point in $S$. 
We fix three distinct points $p_1,p_2$, and $p_3$ on the boundary of $\Omega$ so that when we travel counterclockwise along $\dee \Omega$, the indices are encountered in the increasing order.

Let $\{\dee B_i\}_{i=1}^{\infty}$ be an indexed collection of peripheral circles of $S$.
For each $n=1,2,\dots$, we consider a relative circle domain 
$$
A_n=\Omega\setminus\bigcup_{i=1}^n\overline{B_i}.
$$
Let $g_n$ be the conformal map from $A_n$ onto a relative circle domain $A_n'$ in the unit disc $\mathbb U^2$, so that under the continuous extension of $g_n$ to the boundary, $\dee \Omega$ corresponds to  $\dee \mathbb U^2$, and points $p_1,p_2,p_3$ are mapped to the points $1,i,-1$, respectively. Such a map $g_n$ is unique. 

As in the proof of Theorem~\ref{T:Blu}, using Proposition~\ref{P:Lbilip} and the Arzel\`a-Ascoli Theorem it follows that a subsequence $(g_{n_j})$ of $(g_n)$ converges locally uniformly to a locally bi-Lipschitz map $g\: S\to S'$, where $S'$ is a relative Schottky set in $\mathbb U^2$. The bi-Lipschitz constants depend only on $\Omega, S, \delta,\sigma$, and a lower bound for  the distance to $\partial\Omega$.
Again by Proposition~\ref{P:Lbilip}, the first derivatives $g_n',\ n=1,2,\dots$, are uniformly locally Lipschitz in $\Omega$, and by the proof of the same proposition they are uniformly locally bounded in $\Omega$. The constants depend only on $\Omega, S, \delta$, and a lower bound for the distance to $\partial\Omega$.
Thus another application of the Arzel\`a-Ascoli theorem gives that for a further subsequence of $(g_n)$, whose index sequence we still denote by $(n_j)$, we have $(g_{n_{j}}')$ converges locally uniformly to a continuous function $h\: S\to\C$, which has to be locally Lipschitz. The Lipschitz constant depends only on $\Omega, S, \delta$, and the distance to $\partial\Omega$. We conclude that for every $p\in S,\ r\leq {\rm dist}(p,\dee\Omega)$, and every $\epsilon>0$, there exists $J\in\N$ with
$$
|g_{n_j}(q)-g(q)|<\epsilon,\quad |g'_{n_j}(q)-h(q)|<\epsilon,\quad q\in \overline{B(p,r/4)}\bigcap S,\ j\geq J.
$$    
We will prove that 
\begin{equation}\label{E:Lim}
\lim_{q\to p,\, q\in S}\frac{g(q)-g(p)}{q-p}=h(p).
\end{equation}

Let $\epsilon>0$ be arbitrary. We assume that a point $q$ is contained in $B(p,r/8)$ and consider a curve $l_{p,q}$ in ${B(p,r/4)}\bigcap S$ that connects $p$ and $q$ and such that ${\rm length}(l_{p,q})\leq\pi |q-p|$. Such a curve exists by Corollary~\ref{C:Circles}. 
Since $(g_{n_j}')$ converges to $h$ on compact sets in $\Omega$ and by choosing $q$ to be sufficiently close to $p$, we may assume that 
$$
|g'_{n_j}(z)-h(z)|<\epsilon/2,\quad |h(z)-h(p)|<\epsilon/2,\quad z\in l_{p,q},\ j\geq J.
$$
Then
\begin{equation}\notag
\begin{aligned}
\bigg|\frac{g_{n_j}(q)-g_{n_j}(p)}{q-p}-h(p)\bigg|&=\bigg|\frac{\int_{l_{p,q}}(g_{n_j}'(z)-h(p))dz}{q-p} \bigg|\\ &\leq\frac{\int_{l_{p,q}}|g_{n_j}'(z)-h(p)|dz}{|q-p|}\\ &\leq\pi\epsilon,\quad j\geq J.
\end{aligned}
\end{equation}
By taking the limit as $j\to\infty$ we obtain 
$$
\bigg|\frac{g(q)-g(p)}{q-p}-h(p)\bigg|\leq\pi\epsilon
$$
for $q$ sufficiently close to $p$, establishing~(\ref{E:Lim}). Since $g$ is locally bi-Lipschitz, $h(p)\neq0$.
Since $g$ is a local homeomorphism, (\ref{E:Lim}) immediately implies that $(g^{-1})'(\tilde p)=1/h(p)$, where $\tilde p=g(p)$.

Applying the same arguments to $\tilde S$, we obtain a relative Schottky set $\tilde S'$ in $\mathbb U^2$ and a locally bi-Lipschitz homeomorphism $\tilde g\: \tilde S\to\tilde S'$ that is differentiable at every point in $\tilde S$. The composition $\tilde g\circ f\circ g^{-1}$ is a locally quasisymmetric map between relative Schottky sets $S'$ and $\tilde S'$ in $\mathbb U^2$. Since $S$ is assumed to have measure zero, it is immediate that both, $S'$ and $\tilde S'$, have measure zero. By Theorem~\ref{T:Circrig}, the map $\tilde g\circ f\circ g^{-1}$ must be the restriction to $S'$ of a M\"obius transformation. 
This implies that $f$ is locally bi-Lipschitz in $S$.
The chain rule completes the proof of conformality and the continuity of the derivative.
\qed

\medskip

\section{Proof of Theorem~\ref{T:Quantlip}.}\label{S:Quantlip}

\no
We use the same notations as in the proof of Theorem~\ref{T:Diff}.
By Lemma~\ref{L:Qcext}, we may assume that $f$ is the restriction of a homeomorphism $F\colon \overline\Omega\to\overline{\tilde\Omega}$ that is $H$-quasiconformal in $\Omega$.
Assume that in the construction of $\tilde g$ in the proof of Theorem~\ref{T:Diff} we used $\tilde A_n=\tilde\Omega\setminus\bigcup_{i=1}^n\overline{\tilde B_i}$,
and for a conformal map $\tilde g_n$ of $\tilde A_n$ onto a relative circle domain $\tilde A_n'$ in $\mathbb U^2$ we had $\tilde g_n(f(p_1))=1, \tilde g_n(f(p_2))=i$, and $\tilde g_n(f(p_3))=-1$. 

Look at $F_n=\tilde g_n\circ F\circ g_n^{-1}$. This is an $H$-quasiconformal map from $A_n'$ onto $\tilde A_n'$, so that its continuous extension to the boundary fixes $1,i,-1$. Using the group of reflections in peripheral circles as in Section~\ref{S:Qc}, we can extend $F_n$ to a global $H$-quasiconformal map of the sphere $\Sph^2$ that fixes $1,i,-1$. The extension will still be denoted by $F_n$ and it takes $A_n'$ onto $\tilde A_n'$. The family of such maps $F_n$ is compact, and thus, by possibly passing to a subsequence, we may assume that $(F_n)$ converges to an $H$-quasiconformal map $G$ of $\Sph^2$ that fixes $1,i,-1$. It is clear that $G$ maps $S'$ onto $\tilde S'$ and its restriction to $S'$ agrees with $\tilde g\circ f\circ g^{-1}$. Thus $G$ restricted to $S'$ is the identity, i.e., $S'=\tilde S'$ and $f=\tilde g^{-1}\circ g$.  

Now, let $p\in S$ be arbitrary and $d\leq {\rm dist}(p,\partial \Omega), {\rm dist}(f(p),\dee\tilde\Omega)$. According to Proposition~\ref{P:Lbilip}, each $g_n$ is $L_1$-Lipschitz in $\overline{B(p,d/4)\bigcap A_n}$, where $L_1$ depends only on $d$.  
Taking the limit gives that $g$ is $L_1$-Lipschitz in $\overline{B(p,d/4)\bigcap S}$.
By Proposition~\ref{P:Comp}, ${\rm dist}(g_n(p),\dee\mathbb U)\geq\epsilon$, where $\epsilon>0$ depends only on $\Omega, S, \delta$, and $d$.
Applying Proposition~\ref{P:Lbilip} again, we conclude that $\tilde g_n^{-1}$ is $L_2$-Lipschitz in $\overline{B(g_n(p),\epsilon/4)\bigcap A_n'}$, where $L_2$ depends only on $\epsilon$ and $\sigma$. Taking the limit gives that $\tilde g^{-1}$ is $L_2$-Lipschitz in $\overline{B(g(p),\epsilon/4)\bigcap S'}$. Thus $f$ is $L_1L_2$-Lipschitz in $\overline{B(p,\min\{d/4, \epsilon/(4L_1)\})\bigcap S}$. 

Proposition~\ref{P:Lbilip} implies that $\tilde g_n$ is $L_3$-Lipschitz in $\overline{B(f(p),d/4)\bigcap \tilde A_n}$, where $L_3$ depends only on $d$. Thus $\tilde g$ is $L_3$-Lipschitz in $\overline{B(f(p),d/4)\bigcap \tilde S}$. Choosing $n$ large enough, we may assume that ${\rm dist}(\tilde g_n(f(p)), g_n(p))<\epsilon/8$. Therefore, applying Proposition~\ref{P:Lbilip} once again, we have that $g^{-1}_n$ is $L_4$-Lipschitz in $\overline{B(\tilde g_n(f(p)),\epsilon/8)\bigcap \tilde A_n'}$, where $L
_4$ depends only on $\epsilon$ and $\sigma$. Passing to the limit gives that $g^{-1}$ is $L_4$-Lipschitz in $\overline{B(g(p),\epsilon/8)\bigcap S'}$. Thus we conclude that $f^{-1}=g^{-1}\circ \tilde g$ is $L_3L_4$-Lipschitz in $\overline{B(f(p),\min\{d/4, \epsilon/(8L_3)\})\bigcap \tilde S}$. Combining the above conclusions about Lipschitz properties of $f$ and $f^{-1}$, we obtain the desired bi-Lipschitz property.

The last assertion about the Lipschitz property of the derivative can be established along similar lines using the corresponding parts of Proposition~\ref{P:Lbilip}. The details are left to the reader.
\qed

\end{document}